\documentclass[11pt]{article}

\usepackage{graphicx}
\usepackage{amsmath}
\usepackage{amsfonts}
\usepackage{epsfig}
\usepackage{color}

\newcommand{\be}{\begin{equation}}
\newcommand{\ee}{\end{equation}}
\newcommand{\bes}{\begin{equation*}}
\newcommand{\ees}{\end{equation*}}
\newcommand{\beqn}{\begin{eqnarray}}
\newcommand{\eeqn}{\end{eqnarray}}
\newcommand{\beqns}{\begin{eqnarray*}}
\newcommand{\eeqns}{\end{eqnarray*}}
\newcommand{\lkr}{\left(}
\newcommand{\lkv}{\left[}
\newcommand{\rkv}{\right]}
\newcommand{\rkr}{\right)}
\newcommand{\lfi}{\left\{}
\newcommand{\rfi}{\right\}}

\newcommand{\fr}[1]{(\ref{#1})}

\newcommand{\ph}{\varphi}
\newcommand{\del}{\delta}
\newcommand{\Del}{\Delta}

\newcommand{\af}{\alpha}

\newcommand{\ga}{\gamma}
\newcommand{\te}{\theta}
\newcommand{\om}{\omega}
\newcommand{\lam}{\lambda}

\newcommand{\sig}{\sigma}

\newcommand{\Om}{\Omega}

\newcommand{\vs}{\varsigma}
\newcommand{\vro}{\varrho}

\newcommand{\EE}{\ensuremath{{\mathbb E}}}

\newcommand{\II}{\ensuremath{{\mathbb I}}}

\newcommand{\PP}{\ensuremath{{\mathbb P}}}

\newcommand{\RR}{\ensuremath{{\mathbb R}}}

\newcommand{\Var}{\mbox{Var}}

\newcommand{\sign}{\mbox{sign}}

\newtheorem{theorem}{Theorem}
\newtheorem{lemma}{Lemma}
\newtheorem{corollary}{Corollary}

\newtheorem{remark}{Remark}

\newcommand{\bY}{\mathbf{Y}}

\newcommand{\bldx}{\mathbf{x}}

\newcommand{\bte}{\mbox{\mathversion{bold}$\te$}}

\newcommand{\bxi}{\mbox{\mathversion{bold}$\xi$}}

\newcommand{\calM}{\mathcal{M}}

\newcommand{\calP}{\mathcal{P}}

\newcommand{\calN}{{\cal N}}

\newcommand{\hh}{\widehat{h}}
\newcommand{\hj}{\widehat{j}}
\newcommand{\hu}{\widehat{u}}

\newcommand{\hPhi}{\widehat{\Phi}}
\newcommand{\tPhi}{\widetilde{\Phi}}
\newcommand{\feta}{f_{\eta}}
\newcommand{\fzeta}{f_{\zeta}}
\newcommand{\Phimu}{\Phi_{\mu}}
\newcommand{\hPhimu}{\widehat{\Phi}_{\mu}}
\newcommand{\tPhimu}{\widetilde{\Phi}_{\mu}}
\newcommand{\hPhimuh}{\widehat{\Phi}_{\mu,h}}
\newcommand{\hPhimuch}{\widehat{\Phi}_{\mu,\check{h}_n}}

\newcommand{\us}{u^{*}}
\newcommand{\ts}{t^{*}}
\newcommand{\fs}{f^{*}}
\newcommand{\fos}{f_0^{*}}
\newcommand{\phs}{\ph^{*}}
\newcommand{\phis}{\ph^{*}}
\newcommand{\qs}{q^{*}}
\newcommand{\gs}{g^{*}}
\newcommand{\hqs}{\widehat{q^{*}}}
\newcommand{\fetas}{f_{\eta}^{*}}
\newcommand{\fzetas}{f_{\zeta}^{*}}
\newcommand{\psis}{\psi^{*}}
\newcommand{\fo}{f_0}

\newcommand{\psisom}{\psis_{m1}}
\newcommand{\psistm}{\psis_{m2}}

\newcommand{\Cdel}{C_{\del}}
\newcommand{\fdel}{f_{\del}}
\newcommand{\fdels}{f^{*}_{\del}}
\newcommand{\Ks}{K^{*}}

\newcommand{\qdels}{q^{*}_{\del}}

\newcommand{\Cbg}{C_{b g}}
\newcommand{\Cbs}{C_{b s}}

\newcommand{\ujo}{u_{j1}}
\newcommand{\ujt}{u_{j2}}
\newcommand{\hujo}{\widehat{u_{j1}}}
\newcommand{\hujt}{\widehat{u_{j2}}}
\newcommand{\vjo}{v_{j1}}
\newcommand{\vjt}{v_{j2}}
\newcommand{\hvjo}{\widehat{v_{j1}}}
\newcommand{\hvjt}{\widehat{v_{j2}}}
\newcommand{\hv}{\widehat{v}}
\newcommand{\hmu}{\widehat{\mu}}

\newcommand{\Cgo}{C_{g 1}}
\newcommand{\Cgt}{C_{g 2}}
\newcommand{\Cpho}{C_{\ph 1}}
\newcommand{\Cpht}{C_{\ph 2}}
\newcommand{\Caleph}{C_{\aleph}}

\newcommand{\Uhl}{U_{h l}}

\newcommand{\iii}{\int_{-\infty}^{\infty}} 
\newcommand{\ioi}{\int_{0}^{\infty}} 
 
\newcommand{\iof}{\int_{1}^{4}} 
\newcommand{\iov}{\int_{0}^{1}}  

\newcommand{\kn}{{k_n}}  
\newcommand{\mun}{{\mu_n}} 
\newcommand{\fmuk}{f_{\mu, k}}
\newcommand{\qmuk}{q_{\mu, k}}
\newcommand{\Qmuk}{Q_{\mu, k}}
\newcommand{\vo}{\vs_{0}}
\newcommand{\vso}{\vs_{0}}
\newcommand{\vos}{\vs_{0}}

\newcommand{\Ksv}{K_v^{*}}
\newcommand{\Cvro}{C_{\vro}}

\newcommand{\tC}{\tilde{C}}
\newcommand{\tilh}{{\tilde{h}}}

\newcommand{\teo}{{\te_0}}
\newcommand{\ho}{{h_0}}
\newcommand{\jo}{{j_0}}

\newcommand{\hPhhj}{\hPhi_{h_{\hj}}}

\newcommand{\tilj}{\widetilde{j}}
\newcommand{\htj}{h_{\tilj}}
\newcommand{\tj}{\widetilde{j}}

\begin{document}

\title{\bf { Minimax theory of estimation of linear functionals of the deconvolution density 
with or without sparsity }}

\author{{\em   Marianna Pensky}   \\
         Department of Mathematics,
         University of Central Florida  }

\date{}

\bibliographystyle{plain}
\maketitle

\begin{abstract}
The present paper   considers  the   problem of estimating a  linear functional
$\Phi = \iii \ph (x) f(x) dx$ of an unknown deconvolution density $f$ on the basis of $n$ i.i.d. observations, 
$Y_1, \cdots, Y_n$ of $Y = \te + \xi$, where $\xi$ has a known pdf $g$, and $f$ is the pdf of $\te$.
Although various aspects and particular cases of this problem   have been treated by a number of authors, 
the minimax theory for estimation of $\Phi$  has not been developed so far.
In particular, there are no minimax lower bounds for an estimator of $\Phi$  
in the case of an arbitrary function $\ph$.  The general upper bounds for the risk 
cover only the case when  the Fourier transform of  $\ph$ exists,   
and cannot be automatically applied to the case when $\ph$ is not absolutely or square integrable. 
Moreover, no theory exists for estimating $\Phi$ in the case when  the vector of observations $\bte = (\te_1, \cdot, \te_n)$ is sparse.
In addition, until now, the related problem of estimation of functionals  $\Phi_n = n^{-1} \sum_{i=1}^n \ph(\te_i)$ 
in indirect observations has been treated as a separate problem with no connection to estimation of $\Phi$.

 The objective of the present paper is to fill in these gaps,  to develop the  general minimax theory 
of estimation of $\Phi$, and to relate this problem to estimation of $\Phi_n$. 
We offer  a general, Fourier transform based approach to estimation of $\Phi$  (and $\Phi_n$)  and provide  
upper and   minimax lower bounds for the risk in the case  when   function $\ph$ is square integrable. 
Furthermore, using technique of inversion formulas, we  extend the theory to a number of situations 
when  the  Fourier transform of $\ph$ does not exist, but  $\Phi$ can be presented as a functional of the Fourier 
transform of $f$ and its derivatives. The latter enables us to construct minimax estimators 
of the functionals that have never been handled before such as the odd absolute moments and the generalized moments 
of the deconvolution density. 
Finally, we generalize our results to  
the situation when the  vector $\bte$ is sparse and the objective 
is estimating $\Phi$ (or $\Phi_n$) over the nonzero components only. 
As a direct   application of the proposed theory, we automatically   
recover  multiple recent results and obtain a variety of new ones     
such as estimation of the mixing cumulative distribution function, 
 estimation of the  mixing probability density function with classical and Berkson errors and 
 estimation of the  $(2M+1)$-th  absolute moment of the deconvolution density.  
\\

\vspace{2mm} 

{\bf  Keywords and phrases}: { linear functional,  minimax lower bound, sparsity,  deconvolution  }

\vspace{2mm}{\bf AMS (2000) Subject Classification}: {Primary: 62G20.  Secondary:  62G05, 62G07  }
\end{abstract}

\section{Introduction  }
\label{sec:introduction}
\setcounter{equation}{0}

In the present paper we consider  the problem of estimating a  linear functional
\be \label{eq:Phi}
\Phi = \iii \ph (x) f(x) dx 
\ee
of an unknown deconvolution density $f$ on the basis of observations $Y_i$, $i=1, \cdots, n$, where
\be \label{obser}
Y_i = \te_i +   \xi_i, \quad i=1, \cdots, n.
\ee 
Here, $\te_i$ are i.i.d. random variables with   unknown pdf $f$,  
and $\xi_i$ are i.i.d. random errors with a known pdf  $g$. The variable $Y$ in this case is the mixture of
$\te$ and $\xi$, and the density $f$ is sometimes referred to as the {\it mixing} density.
Vector $\bte = (\te_1, \cdots, \te_n)$ in \fr{obser} may be sparse in the sense that, on the average, it has only $k_n$
non-zero elements where $k_n/n \to 0$ as $n \to \infty$.

Note that the problem of estimating $\Phi$ in \fr{eq:Phi}  appears in many contexts. If $\ph (x) = \del (x-x_0)$, then $\Phi$ is the 
value of the unknown deconvolution density $f$ at the point $x_0$, estimation of which has been studied extensively by 
Butucea and Comte (2009). If $\ph(x) = \II (x< x_0)$, where $\II(\Om)$ denotes the indicator of a set $\Om$,
then problem \fr{eq:Phi} reduces to estimation of  the mixing distribution function $\Phi = F(x_0)$ at   $x_0$  examined 
by Dattner {\it et al.} (2011). If $\ph (x) = e^{i \om_0 x}$, then $\Phi = \widehat{f} (\om_0)$, 
the characteristic function of the mixing distribution at $\om = \om_0$. If $\ph (x) = f_{\eta} (x- x_0)$, where $f_{\eta}$  is a known pdf, 
then $\Phi = \Phi(x_0)$ is itself a convolution density at a point $x_0$, as considered by Delaigle (2007).
Finally, if $\ph(x) = x^k$ or $\ph(x) = |x|^{2M+1}$, then $\Phi$ is, respectively, the $k$-th moment or the $(2M+1)$-th absolute moment of
the mixing density $f$.

In addition, the problem of estimating $\Phi$ in  \fr{eq:Phi} can be related to 
estimating  of functionals in indirect observations 
\be \label{eq:Phin}
\Phi_n  = \frac{1}{n} \sum_{i=1}^n \ph(\te_i).
\ee 
A particular case of this problem, with $\ph(x) = |x|$ and $\xi_i$ being i.i.d. standard Gaussian errors,
 has been  recently considered by  Cai and Low (2011).   
Problems of the form \fr{eq:Phin} appear whenever one is interested in estimation
of some characteristic which cannot be observed directly. For example, it is well known in nutritional research
(see, e.g.,  
Blundell  (1998)),   
that   people routinely miss-report their  caloric intake 
as well as quantities of various food items which they consume.  In order to accurately estimate 
the percentage of calories coming from, for example, fruit and vegetables, or the probability that a person
consumes 25\% more calories than the recommended amount, one needs to take these random errors 
into account. Similar situations occur in astronomy, where all measurements are subject to high levels of
instrumental noise (see, e.g. Jaffe  (2010)).

Indeed, if $\te_i$, $i=1, \cdots, n$ are i.i.d. with the pdf $f$ and  
 $\EE |\ph(\te)|  < \infty$, then $\Phi_n$ in \fr{eq:Phin}
can be viewed as an ``estimator'' of $\Phi = \EE \ph(\te)$ in \fr{eq:Phi} on the basis of
``observations'' $\te_i$. 
Moreover,  as long as $\EE |\ph(\te)|^2  < \infty$, one has $\EE (\Phi_n - \Phi)^2 \leq n^{-1} \EE |\ph(\te)|^2$. Hence,
the minimax  risks for estimating   $\Phi_n$ and $\Phi$ are equivalent up to the $C n^{-1}$
additive term. Therefore, the upper bounds and the minimax lower bounds for the risks of both
estimators will coincide up to, at most, a constant factor.


In the case when vector $\bte$ is sparse and has, on the average, only $k_n$ non-zero components, observations \fr{obser}  can be viewed as
heterogeneous sparse mixture, which was studied  by a number of authors (see, e.g., 
Donoho and Jin (2004), Cai {\it et al. } (2007),  Hall and Jin  (2010)  among others). 
Models of this type appear in various areas of applications such as 
 early detection of bio-weapons use, detection of covert communications,
 meta-analysis with heterogeneity (see, e.g., Donoho and Jin (2004)),
or  detecting stellar occultations by Kuiper Belt objects  (see, e.g., Liang {\it et  al.}  (2004)).
In particular, in the case when $\xi_i$ are i.i.d. Gaussian errors, research has been focused on  
testing multiple hypotheses that $\te_i =0$, and estimating non-zero components $\te_i$, $i=1, \cdots, n$.
 
If the signal is present ($k_n >0$) and $k_n$ is known, the problem of interest is to estimate 
some characteristics of nonzero elements of $\bte$. The latter problem can be summarized as 
the problem of estimating 
\be \label{Phisp}
\Phi_{\kn}  = \frac{1}{\kn} \sum_{i=1}^n \ph(\te_i)\ \II(\te_i \neq 0).
\ee 
Cai and Low (2011) considered estimation of \fr{Phisp} when $\ph(x)=|x|$, the errors are Gaussian 
and $k_n = n^{\nu}$, $0<\nu <1$. They  concluded that consistent estimation is impossible if $\nu \leq 1/2$, while 
estimation yields the same minimax rates as in the non-sparse case for $\nu > 1/2$. The  question of interest  is whether the same will 
happen in general, or whether this phenomenon is due to the type of the functional (the first absolute moment)
or the type  of errors (Gaussian) studied in the paper. 
Note that, in the case of sparse  vector $\bte$, the pdf $f$ of $\te_i$ can be written as 
\bes 
f (x) = \mun \fo (x) + (1 - \mun) \del(x), \quad \mun = \kn/n = n^{\nu -1},
\ees 
where $\fo (\te)$ is pdf of the nonzero entries of $\bte$, and   $\Phi_{\kn}$ in \fr{Phisp} corresponds to 
\be \label{Phi_mu}
\Phi_{\mu} = \iii \ph (x) f_0 (x) dx.
\ee
 

In spite of its great importance, surprisingly, the general  problem of estimation of a  linear functional 
of the deconvolution density has not been thoroughly investigated. 
In the non-sparse case, the problem of   estimation of  linear functional \fr{eq:Phi} of the mixing density 
with a square integrable function  $\ph$ has been addressed by Butucea and Comte (2009) who derived the upper bounds for 
the mean squared risk for a variety of estimation scenarios, and constructed adaptive estimators that attain them 
(up to, at most, a logarithmic factor). However, the minimax lower risk bounds  have been derived only in the case 
when $\ph (x) = \del (x-x_0)$ due to technical difficulties. Hence, although it is intuitively clear that 
the estimators constructed by Butucea and Comte (2009) attain (up to, at most, a logarithmic factor) 
the minimax lower bounds for the risk for other choices of $\ph (x)$, to the best of our knowledge, this has never been proved. 
Moreover, the general study of  Butucea and Comte (2009) has been confined to the situation 
when Fourier transform of function $\ph(\te)$ exists,
which does not allow one to apply the theory when $\Phi$ is, for example, a mixing  cdf at a point $x_0$. 
Furthermore, as far as we know, there has never been a  study of estimation of a  linear functional
of the deconvolution density in the sparse setting. For example, one would like to know for which values 
of $\nu$ one can estimate functionals \fr{Phisp} and \fr{Phi_mu} consistently and what the optimal convergence rate 
would be. 
The existing   results for estimating of linear functionals in a general inverse problems setting, 
derived by Goldenshluger and Pereverzev (2000) and  Math\'{e} and Pereverzev  (2002) focus   on the upper bounds for the
risk, do not involve the sparse setting  and, in addition,  cannot be easily adapted 
to estimation of functionals  \fr{eq:Phi} or \fr{eq:Phin}.

The purpose of the present paper is to   to fill in the existing gaps and to advance the theory of estimation of linear functionals
of the deconvolution density. In particular, the paper accomplishes  several key goals: 

\begin{enumerate}
\item
 Derivation of minimax lower bound  for the risk of an estimator  of a general linear 
functional of the deconvolution density. These bounds have not been obtained before and they confirm that
the estimators obtained via Fourier transform are indeed minimax optimal.

\item
Estimation of linear functionals \fr{eq:Phi} when 
function $\ph$ is not necessarily integrable or square integrable
using inversion formulas. This is a completely new technique.

\item
Application of those methodologies to estimation of   functionals 
of the form \fr{eq:Phin}, which  allows one to obtain estimators of those functionals
with no additional effort.   This part, although technically very simple,  is
philosophically important.

\item
Estimation of the linear functionals of the form  \fr{Phi_mu}  (or \fr{Phisp}) in
the case when deconvolution density $f$ (or vector $\bte$) is sparse. 
The advantage of the approach in this paper is that, by introducing an effective sample size
and bias correction, we largely reduce this problem to estimation of a linear functional in a 
general set up, thus, significantly advancing the existing theory. 
\end{enumerate}

Below we provide a more detailed overview of the paper.
We start the paper with the study of the case when Fourier transform of function $\ph (x)$ in \fr{eq:Phi}   
exists in a regular sense.  We refer to this situation as the  
{\bf standard case} in comparison with the situations considered in Section \ref{sec:challenge}
where   $\Phi$ is represented as a linear functional Fourier transform $\fs$ of $f$ using some inversion formula.
We complete the theory of Butucea and Comte (2009)  
by providing the exact expressions for the upper bounds for the risk in the case when $f$ belongs to a Sobolev ball 
(Section \ref{sec:2.2}) and establish  the matching  minimax lower bounds for the risk for a general function $\ph (x)$ 
(Section \ref{sec:2.4}). This confirms  that the estimators derived in Butucea and Comte (2009)  are indeed asymptotically optimal 
(or near optimal up to a logarithmic factor). In Section~\ref{sec:2.3}, we also explain how one can carry out 
an adaptive choice of the bandwidth using Lepskii's  method.
%
As an application of our methodology, we consider 
pointwise  estimation of the mixing density with classical and Berkson errors  (Section~\ref{sec:2.5}), thus,
advancing the theory developed by Delaigle (2007).

Next, we expand our approach to incorporate estimation of functionals of the form  \fr{eq:Phi} where function 
$\ph(x)$ does not have the  Fourier transform in a regular sense. In this case, functional $\Phi$ can often  be represented 
via the Fourier transform of the deconvolution  density or its derivatives, using technique of inversion formulas (Section \ref{sec:challenge}). 
In Section~\ref{sec:3.1} we provide a general formulation and some examples of the functionals that can be reduced to this form.
Section~\ref{sec:3.2} deals with the construction of estimators of  linear functionals of the real and the imaginary parts 
of the Fourier transforms of   derivatives of the deconvolution  density. It  also evaluates the risk bounds for those estimators.  
Section~\ref{sec:3.3} provides examples of construction of functionals by using this methodology. In particular, 
we construct the estimators for the odd absolute moments of a deconvolution density, and also of a functional of the form 
$\iii \te^m (\te^2+1)^{-1} f(\te) d\te$ with $m \geq 2$.    Minimax upper and lower bounds for the  risks
of those estimators are also derived. 
Also, as a result of application of the general theory, we immediately obtain an estimator of the mixing cdf 
studied by Dattner {\it  et al.} (2011) and the first absolute moment of the mixing density,
and the as well as the minimax lower and upper bounds for the risks of those estimators.  
The latter example significantly advances the theory developed in Cai and Low (2011) by generalizing
it to the case of   non-Gaussian errors and the  mixing densities of various degrees of smoothness.

Section \ref{sec:sparse} deals with the situation where vector $\bte$ or deconvolution density $f$   is sparse.
 We propose a general   procedure designed for estimating functionals   \fr{Phisp} and \fr{Phi_mu}
in the sparse case for any function $\ph(x)$ and any kind of error density $g$, and 
 evaluate the upper bounds for the risk over   Sobolev classes (Section \ref{sec:4.1}). We    
construct the  matching (up to, at most, a logarithmic factor) minimax lower bounds for the risk  over those classes 
(Section \ref{sec:4.3}). We discover that convergence rates in this case are determined by the ``effective'' sample size 
$n^{-1}   k_n^{2} = n^{2\nu-1}$ if $k_n = n^\nu$. The latter  proves  that 
conclusion of Cai and Low (2011) that consistent estimation   is impossible if $\nu \leq 1/2$
applies not only to their particular case ($\ph(x)=|x|$,   Gaussian errors) but to any functional and any distribution of errors.

Finally, in Section \ref{sec:simulation} we carry out a finite sample simulation study of estimation 
of the first absolute moment of the mixing density. In particular, we compare the  
Fourier transform based estimator developed in the paper to the estimator of Cai and Low (2011).
We show that, although estimators have similar precisions, the advantage of the estimators developed 
in the present paper is that they are adaptive to the choice of parameters and the smoothness of the unknown 
deconvolution density.
The paper ends with a discussion in  Section \ref{sec:discussion} .

Although, in this paper,  we restrict our attention to the case 
when the Fourier transform $\fs(\om)$ of the unknown mixing density $f$ has polynomial decay as $|\om| \to \infty$,  
with some additional effort,  the results in the paper can be extended to the case 
when $\fs(\om)$ has exponential decay. This, however, is a matter for future consideration.

\section{The minimax upper and   lower bounds for the risk: the standard case }
\label{sec:standard}
\setcounter{equation}{0}

\subsection{Notations and assumptions}

For any function $t(x)$, we denote its Fourier transform by 
\bes 
\ts (\om) = \iii e^{i \om x} t(x) dx.
\ees
Denote the pdf of  observation  $Y_i$ by $q(y)$, so that $\qs  (\om) = \fs (\om) \gs (\om)$.
\\

\noindent
We assume that the mixing density  belongs to the Sobolev ball $f \in \Om_s (B)$ where
\be \label{Sobclass} 
\Om_s (B) = \lfi \ts:\ \iii  | \ts (\om)  |^2 (\om^2 +1)^s d\om \leq B^2 \rfi, \quad s \geq 0.
\ee 
\\

\noindent
We  introduce the following standard assumptions on the known functions $g$ and $\ph$.
\\

\noindent 
{\bf A1.\ }  There exist non-negative constants $\Cgo$, $\Cgt$, $\af$, $\beta$ and $\ga$ 
such that
\beqn
|\gs (\om)| & \geq & \Cgo  (\om^2 +1)^{-\af/2} \exp(- \ga |\om|^{\beta}), \label{glow}\\
|\gs (\om)| & \leq & \Cgt  (\om^2 +1)^{-\af/2} \exp(- \ga |\om|^{\beta}), \label{gup}
\eeqn
where $\af >0$  and $\beta =0$ whenever $\ga =0$.
\\

\noindent 
{\bf A2.\ }  There exist non-negative constants $\Cpho$, $\Cpht$, $a$, $b$ and $d$ 
such that
\beqn
|\phs (\om)| & \geq & \Cpho  (\om^2 +1)^{-a/2} \exp(- d |\om|^{b}), \label{philow}\\
|\phs (\om)| & \leq & \Cpht  (\om^2 +1)^{-a/2} \exp(- d |\om|^{b}), \label{phiup}
\eeqn
where   $b =0$ whenever $d =0$.
\\

\noindent
In what follows, we use the symbol $C$ for a generic positive constant, 
which   takes different values at different places and is independent of $n$.
Also, for any positive functions $a(n)$ and $b(n)$, we write $a(n) \asymp b(n)$ if the ratio 
$a(n)/b(n)$ is bounded above and below by finite positive  constants independent of $n$.
\\

\subsection{Estimation and the upper bounds for the risk}
\label{sec:2.2}

In Section \ref{sec:standard}, we assume that the functional $\Phi$ in \fr{eq:Phi} can be represented as
\be \label{Phiform}
\Phi = \frac{1}{2\pi}\ \iii \fs (\om)\ \phs (-\om) d\om = 
 \frac{1}{2\pi}\ \iii \ \frac{\qs (\om)}{\gs(\om)}\ \phs (-\om) d\om 
\ee
where the integral is absolutely convergent. This happens if, for example, 
$a> 1$ in \fr{philow}, so that $|\phs|$ and $|\ph|$ are  square integrable, however, 
this is true for a wider variety of functions $\ph$ (e.g., $\ph (x) = \del(x-x_0)$
considered in Butucea and Comte (2009)). We refer to this situation as the  
{\bf standard case} in comparison with the situations considered in Section \ref{sec:challenge}
where   $\Phi$ cannot be represented in the form \fr{Phiform}.

Following Butucea and Comte (2009), we estimate  $\Phi$ in \fr{Phiform}
by
\be  \label{Phiest}
\hPhi_h =  \frac{1}{2\pi}\ \iii \ \frac{\hqs (\om)}{\gs(\om)}\  \phs (-\om)\ \II(|\om| \leq h^{-1})\, d\om 
\ee
where
\be \label{hqs}
\hqs (\om) = n^{-1}\ \sum_{j=1}^n e^{i \om Y_j} 
\ee 
is the unbiased estimator of $\qs (\om)$ and $h=0$ if function $|\phs (\om)|/|\gs(\om)|$ has finite $L^2$-norm.
In particular, the   upper bound for the risk of the estimator $\hPhi_h$ 
over the Sobolev class $\Om_s (B)$ 
\bes 
R_n (\hPhi_h, \Om_s(B)) = \sup_{f \in \Om_s (B)} \EE(\hPhi_h - \Phi)^2
\ees
is given by the following inequality
\be \label{upper_risk}
R_n (\hPhi_h, \Om_s (B)) \leq  \frac{\| g \|_{\infty} }{2 \pi n}\ 
\iii \frac{|\phs(\om)|^2}{|\gs (\om)|^{2}} \, \II(|\om| \leq h^{-1}) d\om +
\frac{B^2}{4 \pi^2} \iii \frac{|\phs(\om)|^2}{(\om^2 +1)^{s}} \, \II(|\om|>h^{-1}) d\om
\ee 
where $\| g \|_{\infty} = \sup_x g(x)$.
In order to analyze the right-hand side of \fr{upper_risk} and similar expressions, the following statement is helpful.


\begin{lemma} \label{lem:rates}
Let  $H(n,h) \equiv  H(n, h; A_1, A_2, b,d,\beta, \ga)$ where
\be \label{Hnh}
H(n,h)    = h^{2 A_1} \exp \lkr -2 d h^{-b} \rkr + n^{-1} \int_0^{1/h}  (\om^2 +1)^{A_2} 
\exp \lkr 2 \ga \om^\beta - 2 d \om^b \rkr d\om.
\ee 
Denote
$$
 \tilh_n = \arg\min_{h} H(n,h), \quad \Del_n \equiv \Del_n (A_1, A_2, b,d,\beta, \ga)  = H(n, \tilh_n; A_1, A_2, b,d,\beta, \ga) 
$$
Then
\be   \label{bounds}
\begin{array}{ll}
\Del_n \asymp   n^{-1}, \quad \tilh_n = 0 & \mbox{if} \quad b>\beta,  \\
\Del_n \asymp   n^{-1}, \quad \tilh_n = 0 & \mbox{if}   \quad b=\beta,   d> \ga >0 \\
\Del_n \asymp   n^{-1}, \quad \tilh_n = 0 & \mbox{if} \quad b= \beta, d = \gamma,  A_2 < - 1/2 \\
\Del_n \asymp   n^{-1} \log \log n, \quad \tilh_n = [\log n/(2\ga)]^{-1/\beta} 
& \mbox{if} \quad b= \beta>0, d = \gamma>0,    A_2 = - 1/2 \\
\Del_n \asymp   n^{-1} (\log n)^{\frac{2A_2 + 1}{b}}, \quad \tilh_n = [\log n/(2\ga)]^{-1/\beta}  
& \mbox{if} \quad b= \beta>0, d = \gamma>0,   A_2 > - 1/2 \\
\Del_n \asymp   n^{-1}  \log n, \quad \tilh_n = n^{-\frac{1}{2 A_1}}, 
& \mbox{if} \quad b= \beta=0, d = \gamma=0,   A_2 = - 1/2 \\
\Del_n \asymp   n^{- \frac{2A_1}{2 A_1 + 2 A_2 +1}}, \quad  \tilh_n = n^{-\frac{1}{2 A_1 + 2 A_2 +1}} 
& \mbox{if} \quad  b= \beta=0, d = \ga = 0,\   A_2 > - 1/2 \\
\Del_n \asymp   (\log n)^{- V_1}\  n^{-d/\ga}, \quad \tilh_n = h^* (n)  
& \mbox{if} \quad  b = \beta >0,\   \gamma> d > 0  \\
\Del_n \asymp   (\log n)^{-\frac{2 A_1}{\beta}}\ \exp \lkr - 2d  [h^* (n)]^{-b} \rkr,
\  \tilh_n = h^* (n)  
& \mbox{if} \quad  \beta > b >0,\ d>0,  \gamma>0 \\
\Del_n \asymp   (\log n)^{- \frac{2A_1}{\beta}},   \quad \tilh_n = h^* (n) 
& \mbox{if} \quad  b= d =  0,\ \beta >0, \ga > 0  
 \end{array} 
\ee 
where 
\bes
h^* (n)  = \lkv \frac{1}{2 \ga} \lkr \log n - \frac{2A_1+2 A_2 +1 -\beta}{\beta} \log \log n \rkr \rkv^{-\frac{1}{\beta}},
\quad 
V_1 = \frac{2A_1}{\beta} - \frac{d (2A_1+2 A_2 +1 -\beta)}{\beta\ \ga}.
\ees
\end{lemma}


Validity of Lemma \ref{lem:rates} can be verified by straightforward calculations.
The  upper bounds for the minimax risk follow directly from \fr{bounds}.

\begin{theorem} \label{th:up_bounds}
If $g$ is bounded above, then, under Assumptions A1  and A2 (inequalities \fr{glow} and \fr{phiup} only), 
one derives 
 \be \label{up_bound}
R_n (\hPhi_{\tilh_n}, \Om_s (B)) \leq C \Del_n (s+a+(b-1)/2,  \af -a, b,d,\beta, \ga)
\ee
where the expressions for ${\tilh_n}$ and $\Del_n (A_1, A_2, b,d,\beta, \ga)$ are given by \fr{bounds}.
\end{theorem} 
 
\medskip

\noindent
Theorem \ref{th:up_bounds} allows to derive upper bounds for the risk of $\hPhi_{\tilh_n}$ as an estimator of functional $\Phi_n$ 
defined in \fr{eq:Phin}.

\begin{corollary} \label{cor:upper_Phin}
Let $\te_i$, $i=1, \cdots, n$, in \fr{obser} be i.i.d. with the pdf $f$. If $\ph(\te)$ is uniformly bounded
$|\ph(\te)| \leq \| \ph\|_\infty <\infty$, then under assumptions of Theorem \ref{th:up_bounds},
one has
\be \label{upper_discrete}
R_n (\hPhi_{\tilh_n}, \Phi_n, \Om_s(B)) = \sup_{f \in \Om_s (B)} \EE(\hPhi_{\tilh_n} - \Phi_n)^2 
\leq 2 R_n (\hPhi_{\tilh_n}, \Om_s (B)) + 2 n^{-1}  \|\ph\|_\infty^2
\leq C R_n (\hPhi_{\tilh_n}, \Om_s (B)), 
\ee 
where $R_n (\hPhi_{\tilh_n}, \Om_s (B))$ is given by expression \fr{up_bound} and $C>0$ is a constant independent of $n$.
\end{corollary}

\subsection{Adaptive estimation}
\label{sec:2.3}

Note that in the expressions for the optimal value of bandwidth $\tilh_n$   in Lemma \ref{lem:rates}, parameters 
$a, \af, b, d, \beta$ and $\ga$ are known, the only unknown parameters are $s$ and $B$.
Hence, the only cases when one needs an adaptive choice of bandwidth are the cases where 
$\tilh_n$ depends on $A_1$. However, in every situation except 
the case when  $d=b=\ga=\beta=0$ and $a<\af + 1/2$, one can easily find  an alternative value 
$\hh_n$ of $h$ that delivers a nearly optimal (up to at most a logarithmic factor) convergence rates.

In the case when  $d=b=\ga=\beta=0$ and $a<\af + 1/2$,  one can use the Lepskii method
for construction of $\hh_n$  (see, e.g.,  Lepski  (1991), Lepski  {\it et al.}  (1997)).   
In order to apply the method, consider the set of bandwidths 
\be \label{hjvalues}
{\cal H} = \lfi h_j = n^{-\frac{1}{2\af}} 2^{-j},\ j=0,1, \cdots, J \rfi  \quad \mbox{with} \quad
2^J \leq (\log n)^{-1} n^{\frac{2\af-1}{2\af}},  
\ee
where $J$ is the largest positive integer satisfying inequality above.
Denote
\be \label{hj}
\hj = \min \lfi j: \ 0 \leq j \leq J;\ |\hPhi_{h_j} - \hPhi_{h_k}| \leq C_\Phi n^{-1/2} \sqrt{\log n}\, h_k^{-(\af - a +1/2)},
\quad \forall k = j,  \cdots, J \rfi,
\ee
where $C_\Phi$ is such that 
\be \label{CPhi}
C_\Phi \geq 4 \max \lfi \frac{\Cpht}{\pi \sqrt{\log n}}; \ \frac{\Cpht}{\Cgo} \lkr \frac{16}{3 \pi} + 
\frac{2 \sqrt{\|g\|_\infty}}{\sqrt{\pi}}\rkr \rfi 
\ee
and constants $\Cpht$ and $\Cgo$ appear in \fr{phiup} and \fr{glow}, respectively.
The following statement provides the minimax upper bounds for the risk when
the bandwidth $h_n$ is chosen adaptively, without the knowledge of parameters $B$ and $s$.


\begin{theorem} \label{th:adapt_bounds}
If $g$ is bounded above, then, under Assumptions A1  and A2 (inequalities \fr{glow} and \fr{phiup} only), 
one obtains the following expressions for $\widehat{R_n} \equiv R_n (\hPhi_{\hh}, \Om_s (B))$
\be \label{adapt_bound}
\begin{array}{ll}
\widehat{R_n} \asymp  n^{-1}, \  \hh_n =0 
& \mbox{if} \quad b>\beta,  \\
\widehat{R_n} \asymp  n^{-1}, \  \hh_n =0 
& \mbox{if}   \quad b=\beta,   d> \ga >0 \\
\widehat{R_n} \asymp  n^{-1}, \  \hh_n =0 
& \mbox{if} \quad b= \beta, d = \gamma,  a > \af + 1/2 \\
\widehat{R_n} \asymp  n^{-1} \log \log n, \  \hh_n = \hh_n^* 
& \mbox{if} \quad b= \beta>0, d = \gamma>0,   a = \af + 1/2 \\
\widehat{R_n} \asymp  n^{-1} \log n, \  \hh_n = n^{-\frac{1}{2a-1}} 
& \mbox{if} \quad b= \beta=0, d = \gamma=0,   a = \af + 1/2 \\
\widehat{R_n} \asymp  n^{-1} (\log n)^{\frac{2\af -2a+1}{\beta}}, \  \hh_n = \hh_n^* 
& \mbox{if} \quad b= \beta>0, d = \gamma>0,   a < \af + 1/2 \\
\widehat{R_n} \asymp  n^{- \frac{2s + 2a -1}{2s + 2\af}} \, \log n, \  \hh_n = h_{\hj} 
& \mbox{if} \quad   b= \beta=0, d = \ga = 0,\   a <  \af + 1/2 \\
\widehat{R_n} \asymp  (\log n)^{-\frac{U_1}{\beta}}\  n^{-d/\ga}, \   \hh_n = \hh_n^*
&  \mbox{if} \quad  b = \beta >0,\   \gamma> d > 0  \\
\widehat{R_n} \asymp  (\log n)^{-\frac{U_2}{\beta}}\ 
  \exp \lkr - 2d \lkv \frac{\log n}{2 \ga} \rkv^{b/\beta} \rkr, \   
\hh_n = \hh_n^* 
& \mbox{if} \quad  \beta > b >0,\ d>0,  \gamma>0,  \\
\widehat{R_n} \asymp  (\log n)^{- \frac{2s + 2a -1}{\beta}}, \   
\hh_n = \hh_n^{**}
& \mbox{if} \quad   b= d =  0,\ \beta >0, \ga > 0  
\end{array} 
\ee
Here, $\hh_n^* = [\log n/(2\ga)]^{-\frac{1}{\beta}}$,
$\hh_n^{**} = [\log n/(3\ga)]^{-\frac{1}{\beta}}$,
$\hj$ is defined in \fr{hj} with $C_\Phi$ defined in \fr{CPhi}, 
$U_1 = \min(\beta + 2a + 2s  -1, \beta + 2a - 2\af -1)$ and $U_2 = b+ 2a + 2s-1$.
\end{theorem} 
 
\medskip

\noindent
 Careful comparisons between the rates in Theorems \ref{th:up_bounds} and \ref{th:adapt_bounds} 
reveal that the rates coincide up to a constant except for the cases 
when $b= \beta=0, d = \ga = 0,\   a <  \af + 1/2$, or $b = \beta >0,  \gamma> d > 0$, or $\beta > b >0,  d>0, \gamma>0$;
in the latter cases the rates coincide up to a log-factor of $n$. One can also replace $\tilh_n$ by $\hh_n$
in Corollary \ref{cor:upper_Phin} and obtain adaptive minimax upper bounds for the risk 
of the estimator of $\Phi_n$.

\begin{remark} {\bf (Comparison with Butucea and Comte (2009)). }
{\rm
Note that, since we are looking at the case of Sobolev classes only, we 
can examine the rate of convergence of the non-specific  term $v_n$ used in Butucea and Comte (2009)
and obtain exact results. In particular, in terms of rates of convergence,
there are three   parametric regions  ($b>\beta$), or ($b=\beta, d > \ga$), or ($b= \beta, d = \gamma,   a > \af + 1/2$), 
several nearly parametric regions ($b= \beta, d = \gamma,   a \leq \af + 1/2$), 
two different regions of polynomial rates of convergence ($b = \beta,   \gamma> d > 0$ or
$d = \ga = 0, \quad a \leq  \af + 1/2$), a region of logarithmic rates 
($b = d =0, \beta >0,   \gamma  > 0$) and an ``in-between'' region 
($\beta > b > 0, d>0,  \gamma>0$) where $R_n (\hPhi_{\tilh_n}, \Om(B))$ converges to 
zero faster than $(\log n)^{-B_1}$  and slower than $n^{-B_2}$
for any $B_1>0,\\  B_2>0$. 
}
\end{remark}

\subsection{The lower bounds for the risk}
\label{sec:2.4}

Theorem \ref{th:up_bounds} provides upper bounds for the risk
for any combination of parameters $s, a, b, d, \af, \beta$ and $\ga$.
However, to the best of our knowledge, the minimax lower bounds for the risk of a general linear functional 
have not been obtained so far.  Butucea and Comte (2009)  derived 
those lower bounds only in the simple case when $|\phs (\om)| =1$.
Below, we derive the minimax lower bounds for the risk and 
show that, for a wide range of functions $\ph$, the upper bounds 
\fr{up_bound} match the lower bounds up to a constant or a logarithmic factor of $n$. \\

Denote 
\be \label{minimax}
R_n (\Om_s (B)) = \inf_{\tPhi}\ \sup_{f \in \Om_s (B)} \EE(\tPhi - \Phi)^2
\ee 
where $\tPhi$ is any estimator of $\Phi$ based on observations $Y_1, \cdots, Y_n$.
Then, the following theorem is true.
\\

\begin{theorem} \label{th:lower_bound}
Let $g$ be bounded above and such that function $\gs$ is differentiable and  
\be \label{gsderiv}
  \frac{|(\gs)' (\om)|}{|\gs (\om)|} \leq C_g\, (1 + |\om|)^{\tau}, \ \tau \geq 0,  
\quad \mbox{with} \quad \tau =0 \quad \mbox{if}\ \ga =0. 
\ee 
Let there exist   $\om_0 \in (0, \infty)$ such that,  for $|\om| > \om_0$,  function  $\rho (\om) = \arg(\phs (\om))$ 
is twice continuously differentiable 
with $|\rho^{(j)} (\om)| \leq \rho  < \infty$,\ $j=0,1,2$.
Then, under Assumptions A1  and A2 (inequalities \fr{gup} and \fr{philow} only), one derives
\be \label{low_bound}
R_n (\Om_s (B)) \geq \lfi
\begin{array}{ll}
C\, n^{-1} & \mbox{if} \  b>\beta\ \mbox{or}\ b=\beta,   d> \ga >0,  \\
C\, n^{-1} & \mbox{if} \     b= \beta, d = \gamma,  a \geq \af + 1/2,  \\
C\, n^{-1}   & \mbox{if} \   b= \beta>0, d = \gamma>0,   a < \af + 1/2, \ \ \ \  \\
C\, n^{- \frac{2s + 2a -1}{2s + 2\af}} & \mbox{if} \    d = \ga = 0,\   a <  \af + 1/2 \\
C\, (\log n)^{-\frac{U_3}{\beta}}\  n^{-d/\ga} &
 \mbox{if} \   b = \beta,\   \gamma \geq  d > 0  \\
C\, (\log n)^{-\frac{U_4}{\beta}}\ \exp \lkr - 2d \lkv \frac{\log n}{2 \ga} \rkv^{b/\beta} \rkr &
\mbox{if} \   b < \beta,\ d>0,  \gamma>0, \\
C\, (\log n)^{- \frac{2s + 2a -1}{\beta}} & \mbox{if} \    b= d =  0,\ \beta >0, \ga > 0  
\end{array} \right.
\ee
where $U_3 = 2a + 2s(1 - d/\ga) - 2d \af/\ga+ 8 \beta - (U_\tau +1)d/\ga -1$, $U_4 =  2a + 2s + 8b -1$
and $U_\tau = \min(7\beta - 2\tau -1,\ 5 \beta +1)$.
\end{theorem}

\medskip

\noindent
 Observe that for in the cases when $\beta  < b$, or $\beta =b,  d> \ga >0$ or $\beta =b,  d = \ga, a \geq \af +1/2$,
the rates of convergence in  Theorem  \ref{th:up_bounds} are parametric and, hence, cannot be improved.
In addition, the lower and upper bounds coincide up to a constant if  $ d = \ga = 0,\   a <  \af + 1/2$
or $b= d =  0,\ \beta >0, \ga > 0$; otherwise, they coincide up to a logarithmic factor.
\\

Since any estimator $\tPhi_n$ of $\Phi_n$ defined in  \fr{eq:Phin}  can be viewed as an estimator of $\Phi$,
due to inequality
\bes 
\EE (\tPhi_n - \Phi_n)^2 \geq 0.5\, \EE (\tPhi_n - \Phi)^2 - 2 n^{-1} \| \ph \|^2_\infty,
\ees
Theorem \ref{th:lower_bound} immediately provides the lower bounds for the risk of any estimator $\tPhi_n$ of $\Phi_n$
based on $Y_1, \cdots, Y_n$.

\begin{corollary} \label{cor:lower_Phin}
Let $\te_i$, $i=1, \cdots, n$, in \fr{obser} be i.i.d. with pdf $f$. If $\ph(\te)$ is uniformly bounded
$|\ph(\te)| \leq \| \ph\|_\infty <\infty$, then under assumptions of Theorem \ref{th:lower_bound},
for sufficiently large $n$, one has
\be \label{lower_discrete}
R_n (\Phi_n, \Om_s (B)) = \inf_{\tPhi_n}\ \sup_{f \in \Om_s (B)} \EE(\tPhi_n - \Phi_n)^2
\geq C R_n (\Om_s (B)),   
\ee 
where $R_n (\Om_s (B))$ is given by expression \fr{low_bound} and $C>0$ is a constant independent of $n$.
\end{corollary}
 
\medskip

As an  example of application of theory above, we   solve the problem of pointwise   estimation of the mixing 
density with classical and Berkson errors studied by Delaigle (2007).
\\

\subsection{Pointwise  estimation of the deconvolution  density with classical and Berkson errors }
\label{sec:2.5}

  Consider the situation where one is interested in estimating the pdf $\fzeta$ of the random variable 
$\zeta = \te + \eta$ where $\te$ and $\eta$ are independent, the pdf $\feta$ of $\eta$ is known and 
one has measurements $Y_1, \cdots, Y_n$ of   random variable $Y = \te + \xi$  of the form of \fr{obser}
where the pdf $g$ of $\xi$ is known.
The model was originally introduced by Berkson (1950) in the regression context and subsequently studied 
by Delaigle (2007) who obtained the upper bounds for the integrated mean squared risk.
In particular, if the pdf $q(y)$ of $Y$ is $k$ times continuously differentiable
and is such that $q$ and $q^{(k+1)}$ are square integrable and  $q^{(k+1)}$ is bounded,   
Delaigle (2007)  derived
\bes
\EE \| \widehat{\fzeta} - \fzeta \|^2  \leq \lfi 
\begin{array}{ll}
C n^{-(2k)/(2k+1 - 2 \beta)}, & \mbox{if}\quad |\fetas (\om)/\gs(\om)| \asymp |\om|^{-\beta} \ \mbox{for}\quad  |\om| \to \infty,\\
C (\log n)^{-2k/\beta}, & \mbox{if}\quad |\fetas (\om)/\gs(\om)| 
\asymp |\om|^b \exp(\ga |\om|^\beta) \quad \mbox{for}\ |\om| \to \infty,
\end{array} \right. 
\ees
where $\|\cdot \|$ denotes the $L^2$-norm with respect to the Lebesgue measure and 
the constant $C$ depends on the  density $f$ of each  $\theta$.

The theory developed in this paper allows one to construct an estimator of the pdf $\fzeta$ at a point $x_0$
with no additional effort. 
Let, as before, $f$, $g$ and $q$ be the pdfs of $\te$, $\xi$ and $Y$, respectively.
Then $\fzetas = \fs \fetas = \qs \fetas/ \gs$ and 
\bes
\fzeta (x_0) = \frac{1}{2\pi} \ \iii e^{-i x_0 \om}\ \frac{\qs (\om) \fetas (\om)}  
{\gs(\om)}\, d\om,
\ees
Therefore,  $\phs (\om) = e^{i x_0 \om} \fetas (-\om)$,   so that $|\phs (\om)| = |\fetas (\om)|$. 
The estimator of $\fzeta (x_0)$ is of the  form \fr{Phiest} and Theorems \ref{th:up_bounds} and 
\ref{th:lower_bound} give the upper and the  minimax lower bounds for the  risk of estimating 
$\fzeta$ at a point $x_0$. In addition, Theorem~\ref{th:adapt_bounds}  provides an adaptive estimator 
of $\fzeta (x_0)$ that, to the best of our knowledge, has not been derived so far.

Here, one can observe an  interesting phenomenon   that we obtain a   wider variety of convergence rates here than   
Delaigle (2007) who recovered only parametric, polynomial (with $d=0$) and logarithmic convergence rates.
The   more diverse convergence rates in our case are not due to the fact that we are studying local (pointwise) error
while Delaigle (2007)  was interested in the global one. Indeed, for this particular example, 
with a little effort, our minimax theory can be extended to the situation  of the global error.
The reason for the wider diversity lies in the fact that we impose assumptions on $\gs$ and $\fetas$
separately while Delaigle (2007) considers only the cases when the absolute value of the ratio $|\fetas/\gs|$
grows polynomially or exponentially as $|\om| \to \infty$.


\section{Estimation of   linear functionals by using inversion formulas }
\label{sec:challenge}
\setcounter{equation}{0}
 
\subsection{Formulation and some inversion formulae}
\label{sec:3.1}

Estimation of the linear functional $\Phi$ in \fr{eq:Phi} relies on the fact that $\ph \in L^2(-\infty, \infty)$,
so its Fourier transform exists. It is easy to see that this condition, however, is not necessary for
consistent estimation of $\Phi$. Consider, for   example, estimation of the $m$-th moment of $f(\te)$
\be \label{mth_mom}
\Phi_m = \iii \te^m f(\te) d\te. 
\ee 
Note that if $\psi_m(y)$ is a solution of the equation 
\be \label{psi_eq}
\iii g(y-\te) \psi_m(y) dy = \te^m
\ee 
then $\Phi_m = \iii \psi_m (y) q(y)dy = \EE[\psi_m(Y)]$.
In order to construct $\psi_m (y)$ satisfying equation \fr{psi_eq}, denote 
\be \label{muknuk}
\mu_{k} = \iii \te^{k} g(\te) d\te, \quad \nu_{k} = \iii \te^{k} f(\te) d\te,
\ee
and assume that $\mu_{2m} < \infty$ and $\nu_{2m} < \infty$.
Let $c_m =1$ and $c_k$, $k = 0, \cdots, m-1$, be solutions 
of the  system of linear equations 
$$ 
\sum_{k=j}^m \mu_{k-j}\, c_k =0, \quad j = 0, 1, \cdots, m-1.
$$ 
Then, it is easy to check that
\bes
\psi_m (y) = y^m + \sum_{k=0}^{m-1} c_k\, y^k, 
\ees
and, under assumption that $\mu_{2m} < \infty$ and $\nu_{2m} < \infty$,  $\Phi_m$ can be estimated by 
\be \label{Phimest} 
\hPhi_m = n^{-1} \sum_{l=1}^n \psi_m(Y_l),\quad \mbox{where} \quad \EE \hPhi_m = \Phi_m, \quad
\Var [\hPhi_m] \leq C_m n^{-1} 
\ee 
and  constant $C_m$ depends only on $m$, $\mu_{2m}$ and $\nu_{2m}$.

Note that although we did not use Fourier transform for estimation of $\Phi_m$ and Fourier transform of $x^m$ does not exist 
in a regular sense, it does exist in a sense of generalized functions and is equal to $(-1)^m \del^{(m)} (\om)$ where 
$\del^{(m)} (\om)$ is the $m$-th derivative of the Dirac delta function (see, e.g. Zayed (1996)).  
However, using the Fourier transform  of $\ph$ as a generalized function would  require $f$ to belong to a so called test-function space.
Those spaces are usually very restrictive, like, e.g., commonly used for the Fourier transforms of generalized functions, the 
space of the Schwartz distributions which consists of all infinitely differentiable functions that vanish outside some compact set 
(see, e.g. Zayed (1996)). One, of course, cannot expect the unknown density $f$ to belong to such space and, moreover,
this will make any minimax estimation totally irrelevant. For this reason, instead of using the theory of generalized functions we 
shall use inversion formulas that mimic generalized functions but do not require unreasonable assumptions on the unknown pdf $f$.
Our goal is to represent the functionals of interest as integrals of the Fourier transform of $\fs$ and its derivatives.
Dattner {\it et al.} (2011) used inversion formula of Gil-Pelaez (1951)  for estimation of the cumulative distribution function
at a point, nevertheless, there are many more possible applications of this technique. Below we consider several examples.
\\

{\bf Example 1.\ Pointwise   estimation of the deconvolution  cumulative distribution function. }
  In order to represent the cdf $F(t)$ of $\te$ at a point $t$ we apply 
formula 3.721.1 of Gradshtein and Ryzhik  (1980)
\be \label{ryzik}
\int_0^{\infty} \frac{\sin(ax)}{x}\ dx = \sign(a),
\ee
where $\sign (a)$ is the sign of $a$.
Denote the real and the imaginary part of $z$ by $\Re [z]$ and $\Im [z]$, respectively. 
Then, due to the relation $\II(\te \leq t) = 1/2 - 1/2\, \sign (\te - t)$, the cdf  $F(t)$ can be represented as $F(t) =  0.5 -   0.5\, \Phi(t)$ where  
\beqn 
\Phi(t) & = &    \iii \sign(\te - t) f(\te) d\te = 
 \frac{2}{\pi}   \, \iii \lkv \ioi   \frac{\sin ((\te - t)\om)}{\om} \,  d\om \rkv  f(\te) d\te  \nonumber\\
 & = &  \frac{2}{\pi}   \, \ioi \frac{\cos(t \om)}{\om}\, \lkv \iii \sin (\te \om) f(\te) d\te \rkv d\om  - 
 \frac{2}{\pi}   \, \ioi \frac{\sin(t \om)}{\om}\, \lkv \iii \cos (\te \om) f(\te) d\te \rkv d\om  \nonumber\\
& = &   \frac{2}{\pi}   \, \ioi \frac{\cos(t \om)}{\om}\, \Im  [\fs(\om)] d\om - 
 \frac{2}{\pi}   \, \ioi \frac{\sin(t \om)}{\om}\, \Re [\fs(\om)] d\om.  \label{altcdf}
\eeqn 
\\ 

\medskip

{\bf Example 2. \ Estimation of truncated moments. }
Consider estimation of $\Phi$ of the form \fr{eq:Phi} where $\ph(\te) = \te^m \II (\te >\teo)$ 
or  $\ph(\te) = \te^m \sign (\te-\teo)$. It is easy to see that these two functionals are related as
\bes
\iii \te^m \II (\te >\teo) f(\te) d\te = \frac{1}{2} \Phi_m + \frac{1}{2}\, \iii \te^m \sign (\te-\teo)  f(\te) d\te
\ees
where $\Phi_m$ is the $m$-th moment of $f$ defined in \fr{mth_mom}. Since $\Phi_m$ can be estimated by $\hPhi_m$ defined 
in \fr{Phimest}  with parametric convergence rates,  it is sufficient to consider
estimation of 
\be  \label{Phima}
\Phi (m, \teo)   =   \iii \te^m \sign (\te-\teo)  f(\te) d\te. 
\ee
Using formula \fr{ryzik}, rewrite $\Phi (m, \teo)$ as
\beqns
\Phi (m, \teo) 
& = & \frac{2}{\pi} \ioi  \frac{\cos(\om \teo)}{\om} \iii \sin(\om \te) \te^m f(\te) d\te
- \frac{2}{\pi} \ioi  \frac{\sin(\om \teo)}{\om} \iii \cos(\om \te) \te^m f(\te) d\te.
\eeqns 
Observe that 
\bes 
\iii \sin(\om \te) \te^m f(\te) d\te = i^{-m} \Im \lkv \frac{d^m}{d \om^m} \ \iii e^{i \om \te} f(\te) d\te \rkv 
= \lfi \begin{array}{ll}
(-1)^{k}  \frac{d^{2k}}{d \om^{2k}} \Im [\fs(\om)], &    m = 2k,\\
& \\
(-1)^{k+1}  \frac{d^{2k+1}}{d \om^{2k+1}} \Re [\fs(\om)], &   m = 2k+1, 
\end{array} \right.
\ees 
and, similarly, 
\bes
\iii \cos(\om \te) \te^m f(\te) d\te  
= \lfi \begin{array}{ll}
(-1)^{k}  \frac{d^{2k}}{d \om^{2k}} \Re [\fs(\om)], &    m = 2k,\\
& \\
(-1)^{k}  \frac{d^{2k+1}}{d \om^{2k+1}} \Im [\fs(\om)], &  m = 2k+1. 
\end{array} \right.
\ees
Combining both cases, we obtain 
\be \label{Phima1}
\Phi (m, \teo) = \lfi \begin{array}{ll}
\frac{(-1)^k\, 2}{\pi} \lkv \ioi \frac{\cos(\om \teo)}{\om}\, \frac{d^{2k} \Im [\fs(\om)] }{d \om^{2k}} d\om -
\ioi \frac{\sin(\om \teo)}{\om} \frac{d^{2k+1} \Re [\fs(\om)]}{d \om^{2k+1}}  d \om \rkv,\!&  m = 2k,\\
& \\
\frac{(-1)^{k+1}\, 2}{\pi}  \lkv \ioi \frac{\cos(\om \teo)}{\om}\, \frac{d^{2k+1} \Re [\fs(\om)]}{d \om^{2k+1}}  d\om +
\ioi \frac{\sin(\om \teo)}{\om}\, \frac{d^{2k+1} \Im [\fs(\om)]}{d \om^{2k+1}}  d\om \rkv,\!&  m = 2k+1.
\end{array} \right.
\ee 
\\ 

\medskip

{\bf Example 3.\ Estimation of generalized moments. }
Consider estimation of functionals of the form \fr{eq:Phi} where $\ph(\te) = \te^m u(\te)$
with $u(\te)$ is such that $u(\te) \in L^2 (-\infty, \infty)$ but  $\te^m u(\te) \not\in L^2 (-\infty, \infty)$.
Note that, since   $u(\te) \in L^2 (-\infty, \infty)$ implies, for any $\teo$, that $u(\te) \II(\te >\teo) \in L^2 (-\infty, \infty)$
and $u(\te) \sign(\te-\teo) \in L^2 (-\infty, \infty)$,  we are automatically including all square integrable
discontinuous functions $u(\te)$. On the other hand, since functions $\II(\te > \teo)$ and $\sign(\te-\teo)$ are not square integrable, 
Example 2 is not a particular case of Example 3. 

In order to derive an inversion formula for this functional, denote $f_m (\te) = \te^m f(\te)$ and observe that 
$\fs_m (\om) = i^{-m} \frac{d^m}{d \om^m} [\fs (\om)]$. Therefore,
\be \label{Phiu1}
\Phi_u = \iii \te^m u(\te) f(\te) d\te = \frac{i^{-m}}{2\pi} \ \iii u^* (-\om) \frac{d^m\, \fs (\om)}{d \om^m} \, d\om.
\ee 
 
 
\subsection{Construction of the estimators and evaluation of their risks}
\label{sec:3.2}

Observe that in all three examples above, the linear functionals $\Phi$ can be presented as a combination of two 
 integrals
\be \label{altPhi}
\Phi =   \ioi \psisom (\om) \,   \frac{d^m \Re [\fs(\om)]}{d \om^m}   \, d\om 
+   \ioi \psistm (\om) \,    \frac{d^m \Im[\fs(\om)]}{d \om^m}   \, d\om. 
\ee
where we assume that $f$ satisfies conditions (that, of course, depend on the particular forms of 
$\psisom(\om)$ and $\psistm(\om)$) which guarantee absolute convergence of the integrals in \fr{altPhi}.
In particular, we assume that both $\fs(\om)$ and $\gs(\om)$ are $m$ times differentiable.
Note that regular case corresponds to $m=0$ and $\psisom (\om) = \psistm (\om) = \phs (-\om)/\pi$.

In order to construct an estimator of the functional $\Phi$ in  \fr{altPhi}, we partition the area of integration into 
${\cal A}_1 = [0;1]$ and ${\cal A}_2 = (1, \infty)$ and rewrite $\Phi$ as $\Phi =  \Phi_1 + \Phi_2$ where 
\be \label{Phi12}
\Phi_k =   \int_{{\cal A}_k} \psisom(\om) \, \Re \lkv \frac{d^m \fs(\om)}{d \om^m} \rkv \, d\om 
+  \int_{{\cal A}_k} \psistm(\om) \, \Im \lkv \frac{d^m \fs(\om)}{d \om^m}  \rkv \, d\om, \quad k=1,2. 
 \ee
Without loss of generality, we assume that $g$ is an even function, so that its Fourier transform $\gs(\om)$ is real-valued. 
One can easily extend our study to the case when 
$g$ is an arbitrary pdf; we leave this case for the reader to examine.
Note that using the general Leibniz rule, $\Phi_1$   in \fr{Phi12} can be re-written as 
\be  \label{alt1}
\Phi_1 =  \sum_{j=0}^m {m \choose j}\ \int_0^1  \frac{d^{(m-j)}}{d\om^{(m-j)}}
\lkv \frac{1}{\gs(\om)} \rkv \lkv \psisom(\om) \ujo  (\om)
+ \psistm(\om) \ujt (\om) \rkv\ d\om 
\ee 
where
\be \label{ujdef}
\ujo  (\om) = \Re \lkv \frac{d^j \qs (\om)}{d \om^j} \rkv, \quad
\ujt  (\om) = \Im \lkv \frac{d^j \qs (\om)}{d \om^j} \rkv
\ee 
Denote 
\bes
\vjo  (\om) = \iii y^j q(y) \cos (\om y) dy, \quad
\vjt  (\om) = \iii y^j q(y) \sin(\om y) dy
\ees
and construct their respective unbiased estimators as
\be \label{hv}
\hvjo  (\om) =  n^{-1} \, \sum_{l=1}^n Y_l^j \cos (\om Y_l), \quad
\hvjt  (\om) =  n^{-1} \, \sum_{l=1}^n Y_l^j \sin (\om Y_l).
\ee
By taking derivatives of $\qs(\om)$ under the integral sign, it is easy to check that
the unbiased estimators of $\ujo(\om)$ and $\ujt(\om)$ are, respectively, given by 
\be \label{hu1}  
\hujo (\om) = \lfi \begin{array} {ll}
(-1)^{j/2}\ \hvjo (\om),& \mbox{if $j$ is even} \\
(-1)^{(j+1)/2}\ \hvjt (\om),& \mbox{if $j$ is odd} 
\end{array} \right.
\ee 
\be  \label{hu2}
\hujt (\om) = \lfi \begin{array} {ll}
(-1)^{j/2}\ \hvjt (\om),& \mbox{if $j$ is even} \\
(-1)^{(j-1)/2}\ \hvjo (\om),& \mbox{if $j$ is odd} 
\end{array} \right.
\ee 
Combination of formulae \fr{alt1} -- \fr{hu2} imply that  $\Phi_1$ can be estimated by 
\be \label{hPhi1}
\hPhi_{1} =   \sum_{j=0}^m {m \choose j}\ \int_{0}^1  \frac{d^{(m-j)}}{d\om^{(m-j)}}
\lkv \frac{1}{\gs(\om)} \rkv \lkv \psisom(\om) \hujo  (\om)
+ \psistm(\om) \hujt (\om) \rkv   d\om.
\ee 
In order to estimate $\Phi_2$, using integration by parts, partition $\Phi_2$ in \fr{Phi12} as    
$\Phi_2 =   F_m(1) + \Phi_{20}$,   
where
\beqns 
F_m (1) \! &\! = \!&\!   \sum_{k=1}^m (-1)^k \lkv \frac{d^{k-1} \psisom (\om)}{d\om^{k-1}}\, 
\frac{d^{m-k}  [\Re(\fs(\om)]}{d\om^{m-k}} + \frac{d^{k-1} \psistm (\om)}{d\om^{k-1}}\, 
\frac{d^{m-k}[\Im(\fs(\om)]}{d\om^{m-k}}   \rkv_{\om =1},\\
%
\Phi_{20} \! &\! = \!&\! (-1)^m\,   \int_1^\infty \lkr
\frac{d^m \psisom (\om)}{d \om^m}\  \Re[\fs(\om)] + \frac{d^m \psistm  (\om)}{d \om^m}\  \Im[\fs(\om)] \rkr d\om. 
\eeqns
Again, taking into account that  $\fs(\om) = \qs(\om)/\gs(\om)$ and applying the  general Leibniz rule,  rewrite $F_m(1)$ and $\Phi_{20}$ as  
\beqns
F_m(1)\! &\! = \!&\!   \sum_{k=1}^m \sum_{j=0}^{m-k} (-1)^k {m-k \choose j}  [A_{m,j,k,1}(1) \ujo(1)  + A_{m,j,k,2}(1) \ujt(1)], \\
\Phi_{20} \! &\! = \!&\! (-1)^m\,   \int_1^\infty \lkv \frac{d^{m} \psisom (\om)}{d\om^{m}}\, u_{01}(\om) + 
\frac{d^{m} \psistm (\om)}{d\om^{m}}\, u_{02} (\om) \rkv \, \frac{1}{\gs(\om)}\, d\om,
\eeqns
where $\ujo(\om)$ and $\ujt(\om)$ are defined in \fr{ujdef} and
\be \label{Amjkl}
A_{m,j,k,l}(\om)   =   \frac{d^{k-1} \psis_{lm} (\om)}{d\om^{k-1}} \, \frac{d^{m-k-j}}{d\om^{m-k-j}} \lkr \frac{1}{\gs(\om)} \rkr, \ l=1,2.
\ee
Therefore, we can estimate $\Phi_2$ by $\hPhi_{2h}=  \widehat{F_m}(1)  + \hPhi_{20h}$ where
\beqn \label{hFm1}
\widehat{F_m}(1)  & = &     \sum_{k=1}^m \sum_{j=0}^{m-k} (-1)^k {m-k \choose j}  [A_{m,j,k,1}(1) \hujo(1)  + A_{m,j,k,2}(1) \hujt(1)]
  \\
 \hPhi_{20h} & = &  (-1)^m\, \int_1^{1/h}  \lkv \frac{d^{m} \psisom (\om)}{d\om^{m}}\, \hu_{01}(\om) + 
\frac{d^{m} \psistm (\om)}{d\om^{m}}\, \hu_{02}(\om) \rkv \, \frac{1}{\gs(\om)}\, d\om,\ \ \    \label{hPhi20}
\eeqn
and $\hujo(\om)$ and $\hujt (\om)$ are defined by \fr{hu1} and  \fr{hu2}. 
Finally, we estimate $\Phi$ in \fr{altPhi} by $\hPhi_h = \hPhi_1 +  \widehat{F_m}(1)  + \hPhi_{20h}$ where $\hPhi_1$,
$\widehat{F_m}(1)$ and $\hPhi_{20h}$  are 
evaluated according to \fr{hPhi1}, \fr{hFm1} and \fr{hPhi20}, respectively.

In order to construct an  upper bound for the risk of the estimator $\hPhi_h$, we denote
\be \label{sigj}
\sig_{j1}^2 (\om) =  n\, \Var(\hujo), \quad 
\sig_{j2}^2 (\om) = n\, \Var (\hujt)
\ee
and consider a class of pdfs    
\be  \label{another_class}
\Xi_{s} (B) = \lfi  f:\ \sup_{\om} \lkv\, |\fs(\om)|\, (|\om|^s +1) \,\rkv \leq B_2  \rfi
\ee
Then, the risk of the estimator $\hPhi_h$ is given by the following statement.

\begin{theorem} \label{th:upper_risk_inv}
Assume that functions  $f$ and $g$ are such that 
\be \label{ass_m}
\mu_{2m} < \infty, \quad \nu_{2m} < \infty,
\ee 
where $\mu_k$ and $\nu_k$ are defined in \fr{muknuk}.  
Let also function  $\gs$ be real-valued,  satisfy Assumption A1, be $m$ times   
differentiable and such that, for some $C_{g} >0$ 
\be \label{cond1_momn}
\left| \frac{1}{\gs(\om)}\, \frac{d^j \gs(\om)}{d \om^j} \right| \leq   C_{g}\, (|\om|  + 1)^{j\, \tau}, \ \tau \geq 0, j = 0, \cdots, m, 
\quad \mbox{where} \quad \tau =0 \quad \mbox{if}\quad \ga =0.
\ee
Let $\psisom (\om)$ and $\psistm (\om)$ be such that for some positive $C_{\psi}$ and nonnegative $d, a_m$ and $b$, and 
for $|\om| \geq 1$ 
\be \label{psiup}
\left| \frac{d^j \psis_{mk}(\om)}{d \om^j} \right| \leq C_\psi  (\om^2 +1)^{- a_m (j+1)/2} \exp(- d |\om|^{b}),\ \  
j=0, \cdots, m, \ k=1,2.
\ee
Assume also that there exists an absolute constant $C_\sig$ such that 
\be \label{extra_assump}
\int_0^1 \lkv \frac{d^{m-j}}{d \om^{m-j}} \lkr \frac{1}{\gs(\om)} \rkr \rkv^2\ |\psis_{mk} (\om)|^2\, \sig_{jk}^2 (\om)\, d\om \leq C_\sig,
\ \ j = 0, \cdots, m,\ .
\ee
Let $\hPhi_h = \hPhi_1 +  \widehat{F_m}(1)  + \hPhi_{20h}$ where $\hPhi_1$,
$\widehat{F_m}(1)$ and $\hPhi_{20h}$  are defined in 
 \fr{hPhi1}, \fr{hFm1} and \fr{hPhi20}, respectively.
Then, 
\be \label{risk_inversion}
\EE(\hPhi_h - \Phi)^2 \! \leq \! C \lkv h^{2A} \exp \lkr - 2 d h^{-b} \rkr + n^{-1} \int_1^{h^{-1}} (\om^2 +1)^{\af - (m+1) a_m} 
\exp \lkr - 2 d \om^b + 2 \ga \om^\beta \rkr \, d\om \rkv,\ \ \ \ 
\ee
where $A = (m+1) a_m + s + b-1$ if $f \in \Xi_s(B)$ and $A = (m+1) a_m + s + (b-1)/2$ if $f \in \Om_s(B)$.
Here, $\Om_s(B)$ and $\Xi_s(B)$ are defined by \fr{Sobclass} and \fr{another_class}, respectively.
\end{theorem} 

Depending on the respective values of $b,d,\beta, \gamma, s, m$ and $a_m$, one can obtain convergence rates and optimal bandwidth values $\tilh$
for each combination of parameters using Lemma \ref{lem:rates} in Section \ref{sec:2.2}.
Moreover, application of an equivalent of Theorem \ref{th:adapt_bounds} allows one to obtain an adaptive estimator of $\Phi$.
Note, however, that convergence rates in Theorem~\ref{th:upper_risk_inv}  are not minimax. Indeed, in addition to 
$f \in \Om_s(B)$ or $f \in \Xi_s(B)$, assumption \fr{extra_assump} imposes   additional conditions on $\fs$ that depend, 
in a non-trivial way, on the shape of functions $\psis_{mk}$, $k=1,2$, thus, modifying the class of functions $f$. 
For this reason, one has to derive   upper and   lower bounds for the minimax risk on a case-by-case basis.
We shall consider some examples in the next section.


\subsection{Examples of estimation of linear functionals using inversion formulas }
\label{sec:3.3}

{\bf Example 1 (continuation).\ Pointwise   estimation of the deconvolution  cumulative distribution function. }
Recall that $F(t) =  0.5 -   0.5 \Phi(t)$ where $\Phi(t)$ is defined by formula \fr{altcdf}.
Since $|\sin (x)/x| \leq 1$ for any $x$, one can easily show that both integrals in \fr{altcdf}
are absolutely convergent provided $\iii |\te| f(\te) d\te < \infty$.  Therefore,  
\fr{altcdf} is a particular case of \fr{altPhi} with $m=0$,
$\psisom (\om) = -  2\,\sin(\om)/(\pi \om)$ and $\psistm (\om) =  2\, \cos (\om)/(\pi \om)$,
so that inequality \fr{cond1_momn} is valid. Observe also that \fr{psiup} holds with $m=0$, $a_0 =1$ and $b=d=0$.
Using notations above, one can write 
\bes 
\Phi(t) = \frac{2}{\pi}\, \ioi \frac{\cos(t \om) v_{02} (\om) - \sin(t \om) v_{01} (\om)}
{\om \, \gs(\om)} \, d\om,
\ees
so that 
\bes 
\hPhi_h  (t) = \frac{2}{ \pi} \, \int_0^{h^{-1}}  \frac{\cos(t \om) \widehat{v_{02}} (\om) - \sin(t \om) \widehat{v_{01}} (\om)}
{\om \, \gs(\om)}   \, d\om.
\ees 
Let $g$ satisfy condition  \fr{glow} and $\mu_{2} < \infty$, $\nu_2 < \infty$,  
where $\mu_k$ and $\nu_k$ are defined in \fr{muknuk}. Then, due to   
\bes 
\sig_{01}^2 (\om) = \iii \cos(\om y) q(y) dy \leq 1,
\quad
\sig_{02}^2 = \iii \sin(\om y) q(y) dy \leq \min(1, 2\,\om^2 (\mu_2 + \nu_2)),
\ees
assumption \fr{extra_assump} is satisfied. 
Therefore, one obtains the upper bound \fr{risk_inversion} for the risk 
with $A = s$ if $f \in \Xi_s(B)$ and $A = s + 1/2$ if $f \in \Om_s(B)$. Here, 
$\Om_s(B)$ and $\Xi_s(B)$ are defined by \fr{Sobclass} and \fr{another_class}, respectively. 
Therefore, the upper bounds for the risk as well as the optimal values for bandwidths 
are provided by formula \fr{bounds} with $d=b=0$, $A_1 = A$ and $A_2 = \af  -1$. The proof that the lower bounds are given  by  
formula \fr{low_bound}  with $d=b=0$ and $a = 1$ can be carried out similarly to the proof 
of Theorem \ref{th:lower_bound} in Section \ref{sec:2.3} and Theorem \ref{th:absMmom} below. 
One can check that convergence rates coincide with the ones derived in Dattner {\it et al.} (2011).
Adaptive choice of the bandwidth by the Lepskii method is described in details in Dattner {\it et al.} (2011).
\\ 

\medskip


{\bf Example 3 (continuation).\ Estimation of generalized moments. }
Consider estimation of functional of the form \fr{Phiu1} where $u(\te) = (\te^2 +1)^{-1}$ and $m \geq 2$. 
It is easy to see that, although $u(\te)$ is absolutely and square integrable, function $\te^m\, (\te^2 +1)^{-1}$ is not.
Then, $\us (\om) = \pi e^{-|\om|}$ and, by formula \fr{Phiu1}, obtain
\be  \label{Phiu}
\Phi_u = \iii  \te^m\, (\te^2 +1)^{-1} f(\te) d\te =  i^{-m}  \, \iii e^{-|\om|}\, \frac{d^m\, \fs (\om)}{d \om^m} \, d\om.
\ee   
Hence, formula  \fr{Phiu} is a particular type of \fr{altPhi} with $\psisom (\om) = 2 e^{-\om}$ and $\psistm =0$ if $m$ is even 
and $\psisom (\om) = 0$   and $\psistm =2 e^{-\om}$ if $m$ is odd. 

Note also that, since function $\us (\om)$ does not have a singularity at zero, 
construction of the estimator does not require  partition of $\Phi_u$ into the parts corresponding to $|\om| \leq 1$ 
and $|\om| > 1$. Instead, applying integration by parts directly to the right-hand side of  \fr{Phiu}, obtain
\be    \label{Phiu2} 
\Phi_u =  i^{-m}  \, \lkv \iii  e^{-|\om|}\, [\sign(\om)]^m \fs(\om) d\om - 2 \sum_{j=0}^{\frac{m-2}{2}} (\fs)^{(m-2-2j)} (0) \rkv.
\ee  
Using the   general Leibniz formula for the $l$-th derivative of the product, write $(\fs)^{(l)} (0)$ as  
\bes
(\fs)^{(l)} (0) = \sum_{k=0}^l {l \choose k}\ i^k\,  \mu_k\  \frac{d^{l-k}}{d \om^{l-k}}\lkv \frac{1}{\gs (\om)} \rkv_{\om =0} 
\ees 
where $\mu_k$ is defined in \fr{muknuk}.
Subsequently, estimate $\Phi_u$ by $\hPhi_{u, h} =    {i^{-m}} (\hPhi_{u,1,h} - 2  \hPhi_{u,2})$ where
\beqns 
\hPhi_{u,1,h} & = & \int_{-1/h}^{1/h}  e^{-|\om|}\, [\sign(\om)]^m \, \frac{\hqs(\om)}{\gs(\om)} \, d\om, \\ 
 \hPhi_{u,2} & = & \sum_{j=0}^{\frac{m-2}{2}} \sum_{k=0}^{m - 2 - 2j}{m - 2 - 2j \choose k}    i^k \hmu_k \, 
\frac{d^{m - 2 - 2j -k}}{d \om^{m - 2 - 2j -k}}\lkv \frac{1}{\gs (\om)} \rkv_{\om =0} 
\eeqns
$\hmu_k = \sum_{l=1}^n Y_l^k$ and ${\hqs} (\om)$ is defined in \fr{hqs}.

Let $f$ and $g$ satisfy conditions \fr{ass_m}. Let also $\gs$ satisfy Assumption A1 and inequality \fr{cond1_momn}.
Then $\sig _{jk}^2 (\om)$ are uniformly bounded for $j = 0, \cdots, m$, $k=1,2$,  and inequality \fr{extra_assump}
holds. Moreover, \fr{psiup} is valid with $a_m=0$ and $b=d =1$. Then, Theorem \ref{th:upper_risk_inv} yields
the upper bounds for the risk of the form \fr{risk_inversion}  where $A = s$ if $f \in \Xi_s(B)$ and 
$A = s + 1/2$ if $f \in \Om_s(B)$.  Applying formula \fr{bounds} of   Lemma \ref{lem:rates}  with modifications 
carried out in Theorem \ref{th:adapt_bounds},  obtain the following   upper bounds for the minimax risk 
$\widehat{R_n} \equiv R_n (\hPhi_{u,\hh_n}, \Om_s (B))$ of the adaptive estimator $\hPhi_{u,\hh_n}$:
\bes  
\widehat{R_n}    \leq \lfi
\begin{array}{ll}
C  n^{-1}, \  \hh_n =0 
& \mbox{if} \quad \beta<1 \ \mbox{or}\ \beta =1, \ga <1,  \\
C  n^{-1}\, (\log n)^{2\af +1}, \  \hh_n = (\log n)^{-1},
& \mbox{if} \quad   \beta =1, \ga =1,  \\
C  n^{-\frac{1}{\ga}}\, (\log n)^{-V_1}, \  \hh_n = [\log n/(2\ga)]^{-1},
& \mbox{if} \quad   \beta =1, \ga > 1,  \\
C (\log n)^{-2A}\, \exp \lfi - 2 \lkr \frac{\log n}{2 \ga} \rkr^{\frac{1}{\beta}} \rfi, \ 
\hh_n = [\log n/(2\ga)]^{-1},
& \mbox{if} \quad   \beta >1, \ga >0.   
\end{array}   \right.
\ees
The lower bounds for the minimax risk  can be obtained similarly to  Theorem \ref{th:lower_bound}
and coincide with the minimax upper bounds when  $\beta<1$ or  $\beta =1, \ga <1$, and differ 
from it by a logarithmic factor of $n$ otherwise.

\medskip


{\bf Example 4.\ Estimation of the ${\bf(2M+1)}$-th absolute moment of the deconvolution density. } 
 Consider estimation of a functional of the form 
\be \label{abmom}
\Phi_{2M+1}     = \iii |\te|^{2M+1}  f(\te)d\te.
\ee
Since $|\te|^{2M+1} = \te^{2M+1} \sign(\te)$,
functional \fr{abmom} is a particular case of $\Phi(m, \teo)$ given by \fr{Phima}
with $m = 2M+1$ and $\teo=0$:
\bes   
\Phi_{2M+1}   = (-1)^{M+1} \frac{2}{\pi}   \, \ioi \frac{1}{\om}\ \frac{d^{2M+1} \, \Re [\fs(\om)]}{d \om^{2M+1}} \, d\om, 
\ees
i.e., $\Phi_{2M+1}$ is of the form \fr{altPhi} with $\psis_{2M+1, 1} (\om) = 2\, (-1)^{M+1} (\pi \om)^{-1}$ and $\psis_{2M+1, 2} (\om) = 0$.
Similarly to \fr{Phi12}, re-write $\Phi_{2M+1}$ as $\Phi_{2M+1}=\Phi_{2M+1,1} + \Phi_{2M+1, 2}$ where
$\Phi_{2M+1,1}$ and $\Phi_{2M+1, 2}$ are the portions of $\Phi_{2M+1}$ evaluated over intervals $[0,1]$ and $(1, \infty)$.
Here,
\bes 
\Phi_{2M+1,1} = (-1)^{M+1} \frac{2}{\pi} \, \sum_{j=0}^{2M+1} {2M+1 \choose j} \int_0^1 \frac{1}{\om}\, 
\frac{d^{2M+1-j}}{d\om^{2M+1-j}} \lkv\frac{1}{\gs(\om)} \rkv \ujo (\om) d\om,
\ees
where $\ujo (\om)$ are defined in \fr{ujdef}. Taking into account relations \fr{hu1} and partitioning the sum above 
into the portions with the even and the odd indices, we obtain the following estimator of $\Phi_{2M+1,1}$
\beqn  
\hPhi_{2M+1,1} & = & \frac{2}{\pi} \sum_{k=0}^M  (-1)^{M+k} \, \lfi {2M+1 \choose 2k+1} 
\int_0^1 \frac{\hv_{2k+1,2}(\om)}{\om} \frac{d^{2(M-k)}}{d\om^{2(M-k)}} \lkv \frac{1}{\gs(\om)} \rkv d\om \right. \nonumber\\
& - & \left. {2M+1 \choose 2k} \int_0^1 \frac{\hv_{2k,1}(\om)}{\om} \frac{d^{2(M-k)+1}}{d\om^{2(M-k)+1}} \lkv \frac{1}{\gs(\om)} \rkv d\om \rfi.
\label{Phi2M1es}
\eeqn
Note that for $\sig_{j1}^2 (\om)$    defined in \fr{sigj}, one has
\beqn 
\sig_{2j, 1}^2 (\om) & \asymp & \Var[Y^{2j}\, \cos(\om Y)] \leq \iii y^{4j} q(y) dy, \label{sig11}\\ 
\sig_{2j+1, 2}^2 (\om) & \asymp & \Var[Y^{2j+1}\, \sin(\om Y)]  
\leq \min \lkv \iii y^{4j +2} q(y) dy,\ \om^2\, \iii y^{4j+4} q(y) dy \rkv. \label{sig22}
\eeqn
Hence, condition \fr{extra_assump} is guaranteed by $\mu_{4M+4}<\infty$ and $\nu_{4M+4} < \infty$ where 
$\mu_{k}$ and $\nu_{k}$ are defined in \fr{muknuk}.

Consider now $\hPhi_{2M+1,2}$. Taking into account that, for any $j = 0,1,2, \cdots$, one has 
\bes 
\frac{d^j \psis_{2M+1, 1} (\om)}{d\om^j}  =  \frac{2}{\pi}\ \frac{(-1)^{M+1+j}\, j! }{\om^{j+1}},
\ees
apply \fr{hFm1} and \fr{hPhi20} for this particular case to obtain 
\beqn \nonumber
\hPhi_{2M+1,2,h}& = & \frac{2}{ \pi} \lfi (-1)^M \sum_{k=1}^{2M+1}\ \sum_{j=0}^{2M+1-k} (k-1)! {2M+1-k \choose j} 
\frac{d^{2M+1-k-j}}{d\om^{2M+1-k-j}} \lkv \frac{1}{\gs(\om)} \rkv_{\om=1} \hujo(1)
\right. \\
& - & \left.  \int_1^{1/h}  \hv_{01} (\om)\, \om^{-(2M+2)}\, [\gs(\om)]^{-1}\, d\om   \rfi \label{Phi2M2es} 
\eeqn
Finally, $\Phi_{2M+1}$ can be estimated as $\hPhi_{2M+1,h}=\hPhi_{2M+1,1} + \hPhi_{2M+1, 2,h}$
where $\hPhi_{2M+1,1}$ and $\hPhi_{2M+1, 2,h}$ are given by \fr{Phi2M1es} and \fr{Phi2M2es}, respectively.
Note that condition \fr{psiup} holds with $m = 2M+1$, $a_m =1$ and $b=d=0$.

In particular, if $M=0$, formulae \fr{Phi2M1es} and \fr{Phi2M2es} yield the following estimator for the first absolute moment of $f$ 
\beqn  
\hPhi_{1,1} & = & \frac{2}{\pi}\,  \int_0^1 \lkr \frac{\hv_{01} (\om)\ (\gs)' (\om) }{\om\, (\gs  (\om))^2 } 
+ \frac{\hv_{12} (\om)}{\om\, \gs(\om)} \rkr \, d\om, \label{Phi1est} \\  
\hPhi_{1,2,h} & = &    \frac{2}{\pi}\, \frac{\hv_{01} (1)}{\gs(1)} - 
\frac{2}{\pi}\,  \int_1^{1/h} \frac{1}{\om^2 }\  \frac{\hv_{01} (\om)}{\gs (\om)} \,  d\om.  \label{Phi2est}   
\eeqn  
and $\hPhi_{1,h}=\hPhi_{1,1} + \hPhi_{1, 2,h}$.
\\

In order to derive  upper and  lower bounds for the minimax risk of the estimator $\hPhi_{2M+1,h}$, 
we introduce   the following sets   of functions: 
\be \label{Newclass} 
\Xi_{s} (B_1, B_2) = \lfi  f:\    \iii \te^{4M+4} f(\te) d\te \leq B_1, \quad
\sup_{\om} \lkv\, |\fs(\om)|\, (|\om|^s +1) \,\rkv \leq B_2  \rfi
\ee  
Theorem below provides   upper and  lower bounds for the  minimax  risk of an estimator of 
$\Phi_{2M+1}$ in \fr{abmom} while Corollary \ref{cor:abs-Mmom} produces similar results for the
discrete version of the functional  $\Phi_{2M+1,n} = n^{-1}\, \sum_{i=1}^n |\te_i|^{2M+1}$.

\begin{theorem} \label{th:absMmom}
Let $\mu_{4M+4} \leq B_g  < \infty$ for some positive constant $B_g$ where $\mu_k$ is defined in \fr{muknuk}. 
Let    function  $\gs$  satisfy   condition \fr{cond1_momn}
and also be such that  
\be \label{cond2_Mmom}
\sup_{|\om| \leq 1}\  \left|\frac{1}{\om}\, \frac{d^{2k+1}}{d \om^{2k+1}} \lkr \frac{1}{ \gs (\om)} \rkr \right| 
\leq   C_{g M},  \quad 0 \leq k \leq M.
\ee
If    inequality \fr{glow} in Assumptions A1  holds,  then
\be \label{up_absMmom}
  R_n (\hPhi_{2M+1, \tilh_n}  , \Xi_{s} (B_1, B_2))  
\leq \lfi
\begin{array}{ll}
C\, n^{-1}, \  \tilh_n =0  & \mbox{if} \quad 2 \af - 4 M < 3, \  \beta = \ga = 0 \\
C n^{-1} \log n,\  \tilh_n = n^{-\frac{1}{4M+2}} & \mbox{if} \quad 2 \af - 4 M = 3, \  \beta = \ga = 0 \\ 
C\, n^{- \frac{4M+ 2s + 2}{2s + 2\af-1}}, \  \tilh_n = n^{-\frac{1}{2s + 2\af-1}},   & \mbox{if} \quad  2 \af - 4 M > 3,  \beta = \ga = 0  \\
C\, (\log n)^{- \frac{4M+ 2s + 2}{\beta}}, \   
\hh_n = \tilh_n^{**} & \mbox{if} \quad   \ \beta >0, \ga > 0, 
\end{array} \right.
\ee
where,  
$\tilh_n^{**} =  \lkv \frac{1}{2 \ga} \lkr \log n - \frac{2s + 2\af -1 -\beta}{\beta} \log \log n \rkr \rkv^{-\frac{1}{\beta}}$.

If inequality \fr{gup} in  Assumptions A1 holds, then
\be \label{low_absMmom}
R_n ( \Xi_{s} (B_1, B_2)) 
\geq \lfi
\begin{array}{ll}
C\, n^{-1} & \mbox{if} \quad 2 \af - 4 M \leq  3, \  \beta = \ga = 0  \\
C\, n^{- \frac{4M+ 2s + 2}{2s + 2\af-1}} & \mbox{if} \quad  2 \af - 4 M > 3,  \beta = \ga = 0   \\
C\, (\log n)^{- \frac{4M+ 2s + 2}{\beta}} & \mbox{if} \quad    \beta >0, \ga > 0  
\end{array} \right.
\ee
Here, 
$R_n (\hPhi_{2M+1, \tilh_n}  , \Xi_{s} (B_1, B_2)) = 
\sup_{f \in \Xi_{s} (B_1, B_2)}\, \EE(\hPhi_{2M+1, \tilh_n} - \Phi_{2M+1})^2$
and
$ R_n ( \Xi_{s} (B_1, B_2))  = \\
\inf_{\tPhi}\ \sup_{f \in  \Om_{s} (B_1, B_2)}\, \EE(\tPhi  - \Phi_{2M+1})^2.$
%
\end{theorem} 

Note that the values of $\tilh_n$ are independent of unknown parameters $s$, $B_1$ and $B_2$ if 
$2 \af - 4 M \leq  3$ and $\beta = \ga = 0$. If $\beta >0,  \ga > 0$, one can replace $\tilh_n$ 
by $\hh_n = [\log n/(3 \ga)]^{-1/\beta}$ and the rates of convergence will not change. 
Finally, if $2 \af - 4 M >  3$ and $\beta = \ga = 0$, one can find the value of $\hh_n$ using 
the Lepskii method similarly to the regular case with the price of $\log n$ in the convergence rates.
 
\medskip

\begin{corollary} \label{cor:abs-Mmom}
Let $\te_i$, $i=1, \cdots, n$, in \fr{obser} be i.i.d. with pdf $f$. If $\Phi_{2M+1,n}$ is defined 
by formula \fr{eq:Phin} with $\ph(\te) = |\te|^{2M+1}$, then under assumptions of Theorem \ref{th:absMmom},
for sufficiently large $n$, one has
\begin{align*}
& R_n (\hPhi_{2M+1, \tilh_n}, \Phi_{2M+1,n}, \Xi_s(B_1, B_2))   =   \sup_{f \in \Xi_s (B)} \EE(\hPhi_{2M+1, \tilh_n}  - \Phi_{2M+1,n})^2 
\asymp R_n (\hPhi_{2M+1, \tilh_n}  , \Xi_{s} (B_1, B_2)), \\  
& R_n (\Phi_{2M+1,n}, \Xi_s (B_1, B_2))   =   \inf_{\tPhi_n}\ \sup_{f \in \Xi_s (B)} \EE(\tPhi_n - \Phi_{2M+1,n})^2
\asymp R_n ( \Xi_{s} (B_1, B_2)),
\end{align*}
where $R_n (\hPhi_{2M+1, \tilh_n}  , \Xi_{s} (B_1, B_2))$ and $ R_n ( \Xi_{s} (B_1, B_2))$ are given by \fr{up_absMmom}
and \fr{low_absMmom}, respectively.
\end{corollary}

\medskip

 \begin{remark} \label{class} 
{\bf The choice of the class of functions. }
{\rm 
Observe that we derived the lower and the upper bounds for the risk not for the subset of the Sobolev ball $\Om_s (B_2)$
but rather for the subset of $\Xi_s(B_2)$. This is motivated by our intention to compare our estimator for the first absolute 
moment with the respective estimator of Cai and Low (2011). One can easily obtain    upper and   lower bounds for the minimax risk  
over the set $\Om_{s} (B_1, B_2))$ in a very similar manner.
}
\end{remark}

\begin{remark} \label{cai}
{\bf Relation to Cai and Low (2011).  }
{\rm 
Cai and Low (2011) studied estimation of $\Phi_n$ of the form \fr{eq:Phin} with $\ph (\te) = |\te|$ based on data generated by 
model \fr{obser} where the errors $\xi_i$ are i.i.d. $\calN (0, \sig^2)$ and there are no probabilistic 
assumptions on vector $\bte$. They showed that
\beqn  
\inf_{\tPhi} \sup_{\bte \in \Theta_n (M_0)} \EE(\tPhi_n - \Phi)^2 & \asymp & M_0^2 \lkr \frac{\log \log n}{\log n} \rkr^2,
\label{Cai1} \\
\inf_{\tPhi} \sup_{\bte \in \RR^n} \EE(\tPhi_n - \Phi)^2 & \asymp &  \frac{1}{\log n}, \nonumber 
\eeqn 
where $\Theta_n (M_0) = \lfi \bte:\ \| \bte \|_\infty \leq M_0 \rfi$. 
By employing a state of the art procedure based on Chebyshev and Hermite polynomials, they constructed adaptive estimators that attain 
these convergence rates.  
With the assumption that $\te_i$ are generated independently from pdf $f$, the problem reduces to
estimation of $\Phi_1$ in \fr{abmom}, the first absolute moment of the mixing density. 
Using formulae \fr{Phi1est}  and \fr{Phi2est}, one can construct an estimator $\hPhi_1$ of $\Phi_1$.

Note that, since in the case of Gaussian errors, one has $\af =0$, $\beta =2$ and $\gamma = \sig^2/2$,  
the estimators \fr{Phi1est}  and \fr{Phi2est} are adaptive if $h = \hh_n=  [2 \log n/(3 \sig^2)]^{-1/2}$,   and 
Corollary \ref{cor:abs-Mmom} implies that 
\bes
R_n (\hPhi_{\hh_n}, \Phi_n, \Xi_{s} (B_1, B_2)) \asymp R_n (\Phi_n, \Xi_{s} (B_1, B_2)) 
\asymp (\log n)^{-(s+1)}   
\ees
 where $\Xi_{s} (B_1, B_2) $ is defined in \fr{Newclass}. Since $f$ is a pdf, it is absolutely integrable, so that $s \geq 0$.
Therefore,  convergence rate \fr{Cai1} corresponds to ``the worst case scenario'' where $s=0$ and $f$ is a combination of 
delta functions. The estimator of Cai and Low (2011) addresses this ''worst-case  scenario''  
but is unable to adapt to a more favorable situation where $|\fs (\om)| \to 0$ as $|\om| \to \infty$. 
In addition, our estimator  is more flexible since it is constructed for any type of error density $g$.  
 
On the other hand, the estimator of Cai and Low (2011) does not impose any probabilistic assumptions on $\te_i$,
and, therefore, can be advantageous when the values of $\te_i$, $i=1, \cdots, n$, are not independent.
%
%
}
\end{remark}


\section{The sparse case }
\label{sec:sparse}
\setcounter{equation}{0}
 
\subsection{Estimation procedure and the upper bounds for the risk}
\label{sec:4.1}


The objective of this section  is to estimate the functional $\Phi_{\mu} = \iii \ph (x) f_0 (x) dx$ defined by \fr{Phi_mu}
where $\fo (\te)$ is pdf of the nonzero entries of $\bte$ and 
\be \label{sparseden}
f (x) = \mun \fo (x) + (1 - \mun) \del(x), 
\ee 
where $\mun$ is known. 
Discrete version of \fr{sparseden} appears as the problem of estimation of the functional 
$\Phi_{\kn}$ defined by formula \fr{Phisp}, where  the average  number   $k_n = n^{\nu}$, $0<\nu <1$, of nonzero entries
of vector  $\bte$   in \fr{obser}   is known 
but locations of the zero entries of $\bte$ are not. Observe that $k_n = n^{\nu}$
corresponds to $\mu_n = n^{-1} k_n = n^{\nu -1}$ in \fr{sparseden}.
Again, similarly to the  non-sparse  case, as long as $\EE |\ph(\te)|^2  < \infty$,
one has $\EE (\Phi_{\kn} - \Phimu)^2 \leq \kn^{-1} \EE |\ph(\te)|^2$, so that
the minimax  errors for estimating   $\Phi_\kn$ and $\Phimu$ are equivalent up to the $C \kn^{-1}$
additive term.

Due to \fr{sparseden}, one has $\Phi = \mun \Phimu + (1 - \mun) \ph (0)$, so that the value of $\Phimu$
can be recovered as 
\be \label{Phimuval} 
\Phimu = \mun^{-1}  \Phi   - \mun^{-1} (1 - \mun) \ph (0).
\ee 
Therefore, we estimate  $\Phimu$ by
\be \label{Phimuest}
\hPhimuh =    \frac{\hPhi_h}{\mun}   -   \frac{1 - \mun}{\mun} \lkv \ph (0) - \del_{h} \rkv  \quad \mbox{with} \quad
\del_{h} = \frac{1}{2\pi}\ \iii \phis( -\om) \, \II(|\om|>h^{-1}) d\om,
\ee
where $\hPhi_h$ is defined in \fr{Phiest} and the  correction term  $\del_h$ is a completely known non-random quantity.


In order to justify estimator \fr{Phimuest}, we derive expressions for its variance and   bias.
Since the second term in  \fr{Phimuest} is non-random, $\Var (\hPhimuh) = \mun^{-2} \Var (\hPhi_h)$ 
where $\Var (\hPhi_h)$  is bounded by 
\bes 
\Var (\hPhimuh) \leq   \frac{\| g \|_{\infty} }{2 \pi n \mun^2}\ 
\iii \frac{|\phs(\om)|^2}{|\gs (\om)|^{2}} \  \II(|\om| \leq h^{-1}) d\om. 
\ees 
Since $\fs (\om) = \mun \fos (\om) + (1 - \mun)$, the bias term of $\hPhimuh$ is of the form 
\beqns
\EE \hPhimuh - \Phimu   & = & \frac{1}{2 \pi \mun} \iii \phis(-\om) \fs (\om) \, \II(|\om| \leq h^{-1}) d\om \\
& - & \frac{1 - \mun}{\mun} \lkv \ph (0) - \del_{h} \rkv - \frac{1}{2 \pi} \iii \phis(-\om) \fos (\om)   d\om\\
 & = &   \frac{1}{2 \pi} \iii \phis(-\om) \fos (\om) \II(|\om| > h^{-1})  d\om + \frac{1 - \mun}{2\pi \mun}\ \Delta (h)
\eeqns  
where
\bes
\Delta (h) = 
  \frac{1}{2\pi}\iii \phis(-\om) \, \II(|\om| \leq h^{-1}) d\om - \ph(0) + \del_h =0.
\ees
Therefore,   for any $f_0 \in \Om_s (B)$, where $\Om_s (B)$ is defined in \fr{Sobclass}, one derives:
\bes 
(\EE \hPhimuh - \Phimu)^2 
\leq \frac{B^2}{4 \pi^2} \iii \frac{|\phs(\om)|^2}{(\om^2 +1)^{s}} \, \II(|\om|>h^{-1}) d\om.
\ees
Hence,
\be \label{MSEPhimu}
\EE (\hPhimuh - \Phimu)^2 \leq \frac{\| g \|_{\infty} }{2 \pi n \mun^2}\ 
\iii \frac{|\phs(\om)|^2}{|\gs (\om)|^{2}} \  \II(|\om| \leq h^{-1}) d\om + 
\frac{B^2}{4 \pi^2} \iii \frac{|\phs(\om)|^2}{(\om^2 +1)^{s}} \, \II(|\om|>h^{-1}) d\om.
\ee
Let $n_\mu = n \mun^2 = n^{-1} k_n^2$ be the new, ''effective'' sample size. Then, comparing \fr{MSEPhimu}  
with \fr{upper_risk}, one immediately observes that the upper bounds for the risk of the   estimator $\hPhi_{\mu, h}$ of $\Phimu$
would coincide with the upper bounds for the risk of the estimator $\hPhi_h$ of $\Phi$ in the non-sparse case if the sample size  
$n$ were replaced by the effective sample size $n_\mu$.   Denote  
\be \label{R_n_mun} 
R_{n,\mun} (\hPhimuh, \Om_s(B)) = \sup_{\fo \in \Om_s (B)} \EE (\hPhimuh - \Phimu)^2
\ee 
where $\Om_s(B)$ is defined in \fr{Sobclass}. If $\nu > 1/2$, then $n_{\mu} = n \mu_n^2 = n^{2 \nu -1} >1$, 
so that combination of Theorem \ref{th:up_bounds}
and formula  \fr{MSEPhimu}  immediately yields the upper bounds for the risk.
\\

\begin{theorem} \label{th:up_sparse}
Let $g$ be bounded above and observations be given by model \fr{obser} where 
$f$ is of the form \fr{sparseden} with $\mun = n^{1 - \nu}$,  $\nu > 1/2$.
Then, under Assumptions A1  and A2 (inequalities \fr{glow} and \fr{phiup} only),
one derives
 \be \label{up_sparse}
R_{n,\mun} (\hPhi_{\mu,\check{h}_n}, \Om_s (B)) \leq C \Del_{n^{1-2\nu}} (s+a+(b-1)/2,  \af -a, b,d,\beta, \ga)
\ee
where the expressions for   $\check{h}_n = \tilh_{n^{2\nu-1}}$ and $\Del_{n^{2\nu -1}} (A_1, A_2, b,d,\beta, \ga)$ are given by \fr{bounds}, 
with $n$   replaced by $n^{2\nu -1}$. In addition, one has the following adaptive   minimax 
convergence rates $\widehat{R_{n,\mun}} \equiv R_n (\hPhi_{\hh_n}, \Om_s (B))$
\be \label{up_adapt_sparse}
\begin{array}{ll}
\widehat{R_{n,\mun}} \asymp  n^{-(2\nu -1)}, \  \hh_n =0 
& \mbox{if} \quad b>\beta,  \\
\widehat{R_{n,\mun}} \asymp  n^{-(2\nu -1)}, \  \hh_n =0 
& \mbox{if}   \quad b=\beta,   d> \ga >0 \\
\widehat{R_{n,\mun}} \asymp  n^{-(2\nu -1)}, \  \hh_n =0 
& \mbox{if} \quad b= \beta, d = \gamma,  a > \af + 1/2 \\
\widehat{R_{n,\mun}} \asymp  n^{-(2\nu -1)} \log \log n, \  \hh_n = \hh_n^* 
& \mbox{if} \quad b= \beta>0, d = \gamma>0,   a = \af + 1/2 \\
\widehat{R_{n,\mun}} \asymp  n^{-(2\nu -1)} \log n, \  \hh_n = n^{-\frac{1}{2a-1}} 
& \mbox{if} \quad b= \beta=0, d = \gamma=0,   a = \af + 1/2 \\
\widehat{R_{n,\mun}} \asymp  n^{-(2\nu -1)} (\log n)^{\frac{2\af -2a+1}{\beta}}, \  \hh_n = \hh_n^* 
& \mbox{if} \quad b= \beta>0, d = \gamma>0,   a < \af + 1/2 \\
\widehat{R_{n,\mun}} \asymp  n^{- \frac{(2\nu -1)(2s + 2a -1)}{2s + 2\af}} \, \log n, \  \hh_n = h_{\hj} 
& \mbox{if} \quad   b= \beta=0, d = \ga = 0,\   a <  \af + 1/2 \\
\widehat{R_{n,\mun}} \asymp  (\log n)^{-\frac{U_1}{\beta}}\  n^{-\frac{d(2\nu -1)}{\ga}}, \   \hh_n = \hh_n^*
&  \mbox{if} \quad  b = \beta >0,\   \gamma> d > 0  \\
\widehat{R_{n,\mun}} \asymp  (\log n)^{-\frac{U_2}{\beta}}\ \exp \lkr - 2d \lkv \frac{\log n}{2 \ga} \rkv^{b/\beta} \rkr, \   \hh_n = \hh_n^* 
& \mbox{if} \quad  \beta > b >0,\ d>0,  \gamma>0,  \\
\widehat{R_{n,\mun}} \asymp  (\log n)^{- \frac{2s + 2a -1}{\beta}}, \   
\hh_n = \hh_n^{**}
& \mbox{if} \quad   b= d =  0,\ \beta >0, \ga > 0  
\end{array} 
\ee
Here, $\hh_n^* = [(2\nu -1) \log n/(2\ga)]^{-\frac{1}{\beta}}$,
$\hh_n^{**} = \lkr \log n/4\ga \rkr^{-\frac{1}{\beta}}$,
$\hj$ is defined in \fr{hj} with $C_\Phi$ given by \fr{CPhi} and $n$ replaced by $n^{2\nu-1}$, 
$U_1 = \min(\beta + 2a + 2s  -1, \beta + 2a - 2\af -1)$ and $U_2 = b+ 2a + 2s-1$.
\end{theorem}

\subsection{The lower bounds for the risk}
\label{sec:4.3}

The  upper bounds for the risk in formula \fr{MSEPhimu} and Theorem \ref{th:up_sparse} suggest that, 
for  $0 \leq \nu \leq 1/2$, one has $n_\mu = n \mun^2 = n^{2\nu -1} \leq 1$ 
and construction of a consistent estimator is impossible for any functional of the form \fr{Phi_mu}. 
The next proposition  shows that this, indeed,  is true in a wide variety of situations.
In particular, under mild assumptions, the risk of no estimator can converge to zero faster than $C (n \mun^2)^{-1}$.

\begin{theorem}\label{th:sp-low-param}
Let $f(x)$ be given by \fr{sparseden} and $\Phimu$ be defined by \fr{Phi_mu}.
If, for some $C_I >0$,  there exist two pdfs, $f_1(\te)$ and $f_2 (\te)$ such that
\be \label{ineq1}
I_k = \iii  g^{-1} (x)\  \lkv \iii g(x-\te) f_k(\te) d\te \rkv^2 \,  dx \leq C_I < \infty, \quad k=1,2,
\ee
and 
\be \label{ineq2}
\Del_{1 2} = \iii \ph(\te) [f_1(\te) - f_2(\te)] d\te \neq 0,
\ee
then
\be \label{low_bou1}
\inf_{\tPhimu}\  \sup_{f_1, f_2}\ \EE(\tPhimu - \Phimu)^2 \geq C \min \lkr 1, (n \mun^2)^{-1} \rkr,
\ee
where $\tPhimu$ is any estimator  of $\Phimu$ based on observations $Y_1, \cdots, Y_n$.
\end{theorem}

\noindent
Theorem \ref{th:sp-low-param} implies that, although the proof of the low bounds for the risk in Cai and Low's (2004, 2011) 
depend heavily on the normality assumption, the fact that one needs   $\nu> 1/2$  
in order to consistently estimate $\Phimu$ remains valid, whether $\xi_i$ in \fr{obser}
are normally distributed or not.  Moreover, Theorem \ref{th:sp-low-param}  does not require
function $\ph$ to be integrable or square integrable, so one can apply this theorem easily to a wide  variety of 
functionals. The quantity $n_\mu^{-1}  = (n \mun^2)^{-1}$ acts as the parametric convergence rate 
that cannot be surpassed. Corollary \ref{cor:gauss}  below shows that Theorem \ref{th:sp-low-param} holds in the case 
when $g$ is a Gaussian pdf. A similar calculation can be repeated in the case when, for example,
$g$ is a doubly-exponential pdf.
\\

\begin{corollary}  \label{cor:gauss}
Let $f(x)$ be given by \fr{sparseden} and $\Phi$ be defined by \fr{Phi_mu}.
If $g(x) = \calN (x| 0, \sig^2)$ is a Gaussian pdf   and $\iii |\ph (x)| g(x) dx < \infty$,
then the lower bound \fr{low_bou1} holds. 
 \end{corollary}

\noindent
Corollary \ref{cor:gauss} generalizes the results of Cai and Low   who 
proved the lower bound \fr{low_bou1} in the case of Gaussian errors if $\ph(x) = x$
(Cai and Low (2004)) and $\ph(x) = |x|$ (Cai and Low (2011)).
\\

One would like to derive the lower bounds for the risk of the form \fr{low_bound}
for any combination of function $\ph(\te)$ and $g(\te)$. Unfortunately, we 
succeeded in doing this only in the cases when $g(\te)$ has polynomial descent as $|\te| \to \infty$
(Theorem \ref{th:low_polyn_g}) or when $d=0$ and $\ga >0$ in Assumptions A1 and A2 (Theorem \ref{th:log_risk}).
Denote 
\be   \label{low_sp_def}
R_{n,\mun} (\Om_s (B)) = \inf_{\tPhimu}\ \sup_{\fo \in \Om_s (B)} \EE (\tPhimu - \Phimu)^2
\ee  
where $\tPhimu$ is any estimator of $\Phimu$ based on observations $Y_1, \cdots, Y_n$
and $\Om_s(B)$ is defined in \fr{Sobclass}.
Then, the following theorem is true.
\\

\begin{theorem} \label{th:low_polyn_g}
Let $f(x)$ be given by \fr{sparseden} and $\Phi$ be defined by \fr{Phi_mu}.
Let $g$ be bounded above and such that $|g(\te)| \geq C_{g1} (\te^2 + 1)^{-\vs}$.
Let function $\gs$ be $\vo$ times continuously differentiable, where $\vo$ 
is the closest integer no less than $\vs$,  and  satisfy the following condition
\be \label{gsvsderiv}
  \frac{|d^l  \gs (\om)|}{d \om^l}  \leq C_{g2}\, |\gs (\om)|\, (1 + |\om|)^{\tau l},\ \ 
l=1,2, \cdots, \vo,  \quad \mbox{where} \quad \tau =0 \quad \mbox{if}\ \ga =0. 
\ee 
Let there exist   $ \om_0 \in (0,   \infty)$ such that  function  $\rho (\om) = \arg(\phs (\om))$ is 
$\vo$ times  continuously differentiable for $|\om| \geq \om_0$,
with $|\rho^{(j)} (\om)| \leq \rho  < \infty$,\ $j=0,1,2, \cdots, \vo$.
Then, under Assumptions A1  and A2 (inequalities \fr{gup} and \fr{philow} only), when $\mun = n^{\nu -1}$ 
with $\nu > 1/2$, one derives
\be \label{sparse_low_polyng}
R_{n,\mun} (\Om_s (B))  
\geq \lfi
\begin{array}{ll}
C\, n^{-(2\nu -1)} & \mbox{if} \quad b>\beta,  \\
C\, n^{-(2\nu -1)} & \mbox{if}   \quad b=\beta,   d> \ga >0 \\
C\, n^{-(2\nu -1)} & \mbox{if} \quad b= \beta, d = \gamma,  a \geq \af + 1/2 \\
C\, n^{- \frac{(2\nu -1)(2s + 2a -1)}{2s + 2\af}} & \mbox{if} \quad   d = \ga = 0,\   a <  \af + 1/2 \\
C\, (\log n)^{-\frac{U_5}{\beta}}\  n^{- \frac{d (2\nu -1)}{\ga}} &
 \mbox{if} \quad  b = \beta,\   \gamma> d > 0  \\
C\, (\log n)^{-\frac{U_6}{\beta}}\ \exp \lkr - 2d \lkv \frac{\log n  (2\nu -1)}{2 \ga} \rkv^{b/\beta} \rkr &
\mbox{if} \quad  b < \beta,\ d>0,  \gamma>0 \\
C\, (\log n)^{- \frac{2s + 2a -1}{\beta}} & \mbox{if} \quad   b= d =  0,\ \beta >0, \ga > 0  
\end{array} \right.
\ee
where  
$U_5 = 2a + 2s(1 - d/\ga) - 2d \af/\ga+ 4 \beta(\vo +1) - (U_{\tau, \vo} +1)d/\ga -1$, 
$U_6 = 4 b (\vo +1) + 2a + 2s-1$
and $U_{\tau, \vo} = \min(\beta (4 \vo + 3)  - 2\tau \vo -1,\  \beta (2 \vo + 3) + 2 \vo -1)$.
\end{theorem}

\noindent
Theorem \ref{th:low_polyn_g} confirms that whenever $g$ has polynomial descent, 
convergence rates in Theorem~\ref{th:up_sparse} are indeed optimal within, at most,  a logarithmic factor of the sample size. 
Another general result refers to the case when $d=0$ and $\ga >0$ in Assumptions A1 and A2.

\begin{theorem} \label{th:log_risk}
Let $f(x)$ be given by \fr{sparseden} and $\Phimu$ be defined by \fr{Phi_mu}.
Let $\mun = n^{\nu -1}$ with $\nu > 1/2$ and   inequalities \fr{gup} and \fr{philow}  
in  Assumptions A1  and A2 hold with $d=0$ and $\ga >0$.
Then, 
\be \label{low_s_log} 
R_{n,\mun} (\Om_s (B)) \geq C\, (\log n)^{- \frac{2s + 2a - 1}{\beta}}.  
\ee  
\end{theorem}

\medskip

\begin{remark}{\bf Estimation of functionals of the form \fr{Phisp}. }
{\rm
Note that, for $1/2 < \nu < 1$, one has $\kn^{-1} =  n^{-\nu} <   n^{1-2\nu} = (n \mun^2)^{-1}$. Therefore,  
if $\te_i$, $i=1, \cdots, n$, in \fr{obser} are i.i.d. with the pdf $f$ of the form \fr{sparseden}, for any $h \geq 0$, one has one has 
\beqns 
\EE (\hPhimuh - \Phi_{\kn})^2 &\leq  & 2 \EE (\hPhimuh - \Phimu)^2 + 2  \kn^{-1}  \| \ph \|_\infty^2 \leq 
\EE (\hPhimuh - \Phimu)^2 + 2  n^{-(1 - \nu)}\,  (n \mun^2)^{-1}  \| \ph \|_\infty^2,\\
\EE (\hPhimuh - \Phi_{\kn})^2 &\geq & 1/2\,\EE (\hPhimuh - \Phimu)^2 - 2  \kn^{-1}  \| \ph \|_\infty^2 \geq
1/2\,\EE (\hPhimuh - \Phimu)^2 - 2  n^{-(1 - \nu)}\, (n \mun^2)^{-1}  \| \ph \|_\infty^2.
\eeqns
Hence, if $n$ is large enough, the  upper and the minimax lower bounds for the risks of the estimators of $\Phi_{\kn}$ and $\Phimu$ 
 coincide up to a constant factor.  In particular, the following statement holds.

\begin{corollary} \label{cor:sparse_Phin}
Let $\te_i$, $i=1, \cdots, n$, in \fr{obser} be i.i.d. with pdf $f$ defined in \fr{sparseden}. 
If $\ph(\te)$ is uniformly bounded
$|\ph(\te)| \leq \| \ph\|_\infty <\infty$, then under assumptions of Theorem \ref{th:up_sparse},
one has
\bes 
R_{n,\mun} (\hPhimuch, \Phi_{\kn}, \Om_s(B)) = \sup_{f_0 \in \Om_s (B)} \EE(\hPhimuch - \Phi_{\kn})^2 
\leq C R_{n,\mun} (\hPhimuch, \Om_s(B)), 
\ees 
where $R_{n,\mun} (\hPhimu, \Om_s(B))$,   defined in \fr{R_n_mun}, is  given by expression \fr{up_sparse} 
and $C>0$ is independent of $n$.
Also,   under assumptions of Theorem \ref{th:low_polyn_g} (or Theorem \ref{th:log_risk}) ,
for sufficiently large $n$, one has
\bes 
R_{n,\mun} (\Phi_{\kn}, \Om_s (B)) = \inf_{\tPhi_{\kn}}\ \sup_{f \in \Om_s (B)} \EE(\tPhi_{\kn} - \Phi_{\kn})^2
\geq C R_{n,\mun} (\Om_s (B)),   
\ees 
where $R_{n,\mun} (\Om_s (B))$, defined in \fr{low_sp_def}, is given by expression 
\fr{sparse_low_polyng} (or \fr{low_s_log})  and $C>0$ is independent of $n$.
 \end{corollary}
}
\end{remark}

\section{Simulation study }
\label{sec:simulation}
\setcounter{equation}{0}

 In order to evaluate small sample properties of the estimators presented in the paper we 
carried out a limited simulation study. In particular, we compared our estimator 
$\hPhi_{1,h}=\hPhi_{1,1} + \hPhi_{1, 2,h}$  of the first absolute moment of the mixing density $\Phi_1$, 
where $\hPhi_{1,1}$ and $\hPhi_{1, 2,h}$ are defined in \fr{Phi1est} and \fr{Phi2est}, respectively,  
with the estimator $\hPhi_{CL}$ of Cai and Low (2011) 
which is based on approximation of the absolute value by combination of Chebyshev and Hermite polynomials:
\be   \label{cailow}
\hPhi_{CL} = \sum_{k=1}^{K_*} G_{2k}^* M_0^{1 - 2k} B_{2k}.
\ee   
Here, $M_0$ is such that $|\te_i| \leq M_0$, $i=1, \cdots, n$, $G_{2k}^*$ is the  coefficient for $\te^{2k}$  in the expansion 
of $|\te|$ by Chebyshev polynomials and
\bes
B_{2k} = n^{-1}\, \sum_{i=1}^n H_{2k} (y_i),
\ees
where $H_{2k} (x)$ are Hermite polynomials (with respect to $\exp(-x^2/2)$) of degree $2k$.

We considered three choices for the mixing density $f$:
$f$ is Gaussian    $\calN (m_\te, \sig^2_\te)$, $f$ is uniform on the interval $[a, b]$ and 
$f$ is a  is  combination of delta functions
\be \label{discrete_f}
f (x) = \sum_{k=1}^K b_k \del(x - a_k).
\ee 
Note that, in the latter case, variable  $\te$ is discrete and does not have a pdf in a regular sense.
In particular, we chose $m_\te = 2$ and $\sig_\te = 4$ when  $f$ is Gaussian, 
$a = -2$ and $b=5$ when $f$ is uniform and $K = 3$, $a_1 = -2$, $a_2 =0$, $a_3 =5$ and $b_1 = b_2 = b_3 = 1/3$ 
when $f$ is a  is a combination of delta functions.

For the purpose of comparison, we generated the vector $\bte$ using probability density function (or probability
mass function) $f$  and added a vector of independent  errors $\bxi$ to it to obtain 
the vector of observations $\bY$. Since procedure of Cai and Low (2011) is developed for Gaussian 
errors only, we chose $g$ to be the standard normal pdf  $\calN(0,1)$. We considered four values of sample size,
$n = 500$, $n=200$, $n=100$ and $n = 50$. 
For each of the combinations of vectors $\bte = (\te_1, \cdots, \te_n)$ and $\bY = (Y_1, \cdots, Y_n)$, 
we  evaluated  
\be  \label{absPhi_n}
\Phi_n = n^{-1} \sum_{i=1}^n |\te_i| 
\ee   
and obtained estimators $\hPhi_h$ and $\hPhi_{CL}$ based on observations $\bY$.
While constructing estimator $\hPhi_h$ given by \fr{Phiest}, we use the   value 
$h=\hh_n = [2\, \log n/(3 \sig^2)]^{-1/2}$  with $\sig =1$,
so that estimator $\hPhi_{\hh_n}$  is fully adaptive.    
Estimator $\hPhi_{CL}$ depends on two parameters, $K_*$ and $M_0$. Following Cai and Low (2011), we chose 
$K_* = 0.5\, \log n/\log \log n$. As far as parameter $M_0$ is concerned, we initially chose $M_0$ on the basis 
of our knowledge of $f$. In particular, we used  $M_0= \max  (|a_k|)$ 
when   $f$ is of the form \fr{discrete_f},  $M_0  = \max(|a|, |b|)$ when $f$ is uniform, and $M_0= M_{n,0} = |m_\te| + \sig_\te \sqrt{ 2 \log n}$
if $f$ is Gaussian    $\calN (m_\te, \sig^2_\te)$.  Thus, the estimator $\hPhi_{CL}$ with those choices of $M$ is not adaptive.
Since in practice, the upper bound for the absolute values  $|\te_i|$ is not known, we also studied the case when 
$f$ is of the form \fr{discrete_f} and  as it is suggested in Cai and Low (2011), $M_{n,0}$ is set to $M_{n,0} = \sqrt{2 \log n}$.

 \begin{table} 
\begin{center}
\begin{tabular}{|c| c|c |c| c| c|c|c|c| }
  \multicolumn{9}{ l }{{\sc  Comparison of the  accuracy of the  Fourier transform based  estimator $\hPhi_{\hh_n}$    }}\\
 \multicolumn{9}{ l }{{\sc   and the Cai and Low (2011)   estimator $\hPhi_{CL}$ over 1000 simulation runs   }}     \\
%
 \multicolumn{9}{ l }{  }\\
\hline \hline 
\multicolumn{9}{ |c| }{ $f = [\del(x +2) + \del(x) + \del(x-5)]/3$, $\Phi = 2.33333$,   $M_0=5$ (true value) }\\
\hline
 & \multicolumn{2}{ |c |}{n=500} & \multicolumn{2}{c|}{n=200}  & \multicolumn{2}{c|}{n=100} & \multicolumn{2}{c|}{n=50}\\
\hline
 & $\hPhi_{\hh_n}$ & $\hPhi_{CL}$   & $\hPhi_{\hh_n}$ & $\hPhi_{CL}$  & $\hPhi_{\hh_n}$ & $\hPhi_{CL}$  & $\hPhi_{\hh_n}$ & $\hPhi_{CL}$ \\
 \hline
$\Del $   & 0.02365   & 0.01292   & 0.03727   & 0.02659   &  0.05631  &  0.05225  &  0.09774  &  0.12955 \\   
          & (0.02664) & (0.01698) & (0.05428) & (0.03612) & (0.07795) & (0.07583) & (0.14302) & (0.18895) \\  
\hline 
$\Del_n $ & 0.01606   &  0.00583  & 0.01834   &  0.00816  &  0.02174  &  0.01228 &  0.02938  &  0.06771 \\   
          & (0.01204) & (0.00574) & (0.01862) & (0.01016) & (0.02803) & (0.01681)& (0.03974) & (0.08313) \\  
\hline 
$s(\Phi_n)$ & \multicolumn{2}{ c|}{0.00854 }  &  \multicolumn{2}{c|}{0.02027}  & \multicolumn{2}{c |}{0.04168 }  &  \multicolumn{2}{c|}{0.08186 }  \\  
\hline
\hline
\multicolumn{9}{ l }{  }\\
\multicolumn{9}{ l }{  }\\
\hline 
\hline 
\multicolumn{9}{ |c| }{ $f$ is uniform on $[-2, 5]$, $\Phi = 2.07142$,  $M_0=5$  (true value)   }\\
\hline
 & \multicolumn{2}{c |}{n=500} & \multicolumn{2}{c |}{n=200}  & \multicolumn{2}{c |}{n=100} & \multicolumn{2}{c |}{n=50}\\
\hline
 & $\hPhi_{\hh_n}$ & $\hPhi_{CL}$   & $\hPhi_{\hh_n}$ & $\hPhi_{CL}$  & $\hPhi_{\hh_n}$ & $\hPhi_{CL}$  & $\hPhi_{\hh_n}$ & $\hPhi_{CL}$ \\
 \hline
$\Del $  &  0.00693  &   0.04567  & 0.01566   &  0.06185  &  0.03020  &  0.07095  & 0.06006  &  0.05382 \\   
         & (0.00978) & (0.03314)  & (0.02331) & (0.05857) & (0.04428) & (0.08070) & (0.08486) & (0.07184) \\  
\hline  
$\Del_n $  &  0.00266  &  0.04260  &  0.00582  &  0.04772  &  0.01072  &  0.04927  &  0.02147  & 0.02901 \\   
           & (0.00383) & (0.01804) & (0.00786) & (0.02877) & (0.01533) & (0.04262) & (0.03066) & (0.03884) \\  
\hline
$s(\Phi_n)$ & \multicolumn{2}{ c|}{0.00431 }  &  \multicolumn{2}{c|}{ 0.00992 }  & \multicolumn{2}{c |}{0.01887 }  &  \multicolumn{2}{c|}{ 0.03868 }  \\  
\hline
\hline
\multicolumn{9}{ l }{  }\\
\multicolumn{9}{ l }{  }\\
\hline
\hline
 \multicolumn{9}{ |c| }{ $f$ is Gaussian   $\calN (m_\te, \sig^2_\te)$, $m_\te =2$, $\sig_\te=4$, $\Phi = 3.58237$, 
$M_{n,0} = m_\te + \sig_\te \sqrt{2 \log n}$ (true value) }\\
\hline
 & \multicolumn{2}{c |}{n=500} & \multicolumn{2}{c |}{n=200}  & \multicolumn{2}{c |}{n=100} & \multicolumn{2}{c |}{n=50}\\
\hline
 & $\hPhi_{\hh_n}$ & $\hPhi_{CL}$   & $\hPhi_{\hh_n}$ & $\hPhi_{CL}$  & $\hPhi_{\hh_n}$ & $\hPhi_{CL}$  & $\hPhi_{\hh_n}$ & $\hPhi_{CL}$ \\
 \hline
$\Del $  &  0.01566  & 0.01746   & 0.04382   &  0.03117  &  0.08756  &  0.06384  &  0.16282  & 0.24934 \\   
         & (0.02252) & (0.02324) & (0.06122) & (0.04318) & (0.12636) & (0.09740) & (0.22401) & (0.27030) \\   
\hline
$\Del_n $  &  0.00237  & 0.01030   & 0.00539   &  0.00769  & 0.01104   &  0.01124  &  0.01942  & 0.22841  \\   
           & (0.00338) & (0.00939) & (0.00735) & (0.01047) & (0.01516) & (0.01652) & (0.02673) & (0.16487) \\    
\hline
$s(\Phi_n)$ & \multicolumn{2}{ c|}{ 0.01376 }  &  \multicolumn{2}{c|}{0.03881 }  & \multicolumn{2}{c |}{0.08029 }  &  \multicolumn{2}{c|}{ 0.14125 }  \\  
\hline
\hline
\multicolumn{9}{ l }{  }\\
\multicolumn{9}{ l }{  }\\
\hline
\hline
\multicolumn{9}{ |c| }{ $f = [\del(x +2) + \del(x) + \del(x-5)]/3$, $\Phi = 2.33333$,   $M_{n,0} = \sqrt{2 \log n}$   }\\
\hline
 & \multicolumn{2}{ |c |}{n=500} & \multicolumn{2}{c|}{n=200}  & \multicolumn{2}{c|}{n=100} & \multicolumn{2}{c|}{n=50}\\
\hline
 & $\hPhi_{\hh_n}$ & $\hPhi_{CL}$   & $\hPhi_{\hh_n}$ & $\hPhi_{CL}$  & $\hPhi_{\hh_n}$ & $\hPhi_{CL}$  & $\hPhi_{\hh_n}$ & $\hPhi_{CL}$ \\
 \hline
$\Del $  &  0.02368  &  1.05701  &   0.03770 & 2.40929   &  0.06359  &  4.76195  &  0.09686  &  1.11470 \\   
         & (0.02699) & (0.25291) & (0.04929) & (0.81377) & (0.08691) & (2.20037) & (0.13408) & (1.06931) \\  
\hline 
$\Del_n $  &  0.01590  &  1.05668  &  0.01804  &  2.44297  &  0.02399  &  4.83772  &  0.03042  & 0.91974  \\   
           & (0.01149) & (0.26219) & (0.01812) & (0.91473) & (0.02890) & (2.54537) & (0.04237) & (0.62472) \\  
\hline 
$s(\Phi_n)$ & \multicolumn{2}{ c|}{ 0.00872}  &  \multicolumn{2}{c|}{ 0.01986 }  & \multicolumn{2}{c |}{ 0.04504 }  &  \multicolumn{2}{c|}{ 0.08724 }  \\  
\hline
\hline
\multicolumn{9}{ l }{  }\\
 \end{tabular}
\end{center}
\caption{  The average values (over 1000 runs) of the errors $\Del$ and $\Del_n$ (with the  standard 
deviations of the errors   listed in the parentheses),  and  the average squared  
deviation $s(\Phi_n)$  of $\Phi_n$ from its mean $\Phi$.  
  } 
\label{table1}
\end{table}

We  repeated the process $N=1000$ times.   The accuracies of estimators are determined  
by their average squared deviations $\Del$ and $\Del_n$ from, respectively, $\Phi$ defined in \fr{abmom}
 and $\Phi_n$ defined in \fr{absPhi_n}: 
\bes 
\Del = N^{-1} \sum_{i=1}^N (\hPhi_{\hh_n}(i) - \Phi)^2, \quad
\Del_n = N^{-1} \sum_{i=1}^N (\hPhi_{\hh_n}(i) - \Phi_n(i))^2.
\ees
Here $\hPhi_{\hh_n}(i)$ and $\Phi_n(i)$ are the values of the estimator $\hPhi_{\hh_n}$ and $\Phi_n$ in the $i$-th realization of the simulations.
The standard deviations of $\Del$ and $\Del_n$ are listed below in parentheses.
As a benchmark for the  magnitude of the error, we also report the average squared  
deviation $s(\Phi_n)$  of $\Phi_n (i)$ from its mean $\Phi$:
\bes 
s(\Phi_n) = N^{-1} \sum_{i=1}^N (\Phi_n(i) - \Phi)^2.
\ees 
Results of simulations   are summarized in Table \ref{table1}.

Simulations confirm that, as long as one uses the exact value of $M_0$ in \fr{cailow}, estimators  $\hPhi_{\hh_n}$
and $\hPhi_{CL}$ have very similar precisions: the difference between the average errors  of 
the two estimators is  smaller than respective standard deviations or the average squared deviation 
of $\Phi_n$ from its mean $\Phi$. As it was expected, estimator  
$\hPhi_{CL}$ performs better when $\te$ is a discrete random variable while the set up where $\te$ is a  continuous
random variable benefits $\hPhi_{\hh_n}$. Also, though $\hPhi_{\hh_n}$ is designed to estimate $\Phi$ while $\hPhi_{CL}$
is intended to estimate $\Phi_n$, the error $\Del_n$ turns out to be smaller for $\hPhi_{\hh_n}$ while $\hPhi_{CL}$ 
is somewhat more accurate as an estimator of $\Phi$. The latter confirms that the problems are equivalent
up to a small additive error and that both problems should be treated in a similar fashion.

However, the advantage of the estimator $\hPhi_{\hh_n}$ is that it is adaptive: indeed we fixed the value of  parameter $h$
for all three distributions of $\te$. On the other hand, estimator $\hPhi_{CL}$ is not adaptive since parameter $M_0=M_{n,0}$ 
is chosen on the basis of the pdf $f$ which is not known. Our simulations show that estimator $\Phi_{CL}$ is very sensitive to the choice of $M_0$.
Indeed, when, instead of the true value $M_0$,  we used $M_0 = M_{n,0} = \sqrt{2 \log n}$ suggested in Cai and Low (2011), 
precision of the estimator deteriorated so much   that $\hPhi_{CL}$ stopped being consistent.  
The latter demonstrated the advantage of the estimator $\hPhi_{\hh_n}$ in comparison with $\hPhi_{CL}$. 
Hence, although the estimator $\Phi_{CL}$ is adequate in asymptotic setting, it runs into serious problems for 
moderate values of $n$ that are more likely to occur in practical situations.


\section{Discussion }
\label{sec:discussion}
\setcounter{equation}{0}

This paper significantly advances  the theory of estimation of linear functionals of the unknown deconvolution density. 
In particular,  in the case when the Fourier transform of  $\ph$ exists, we derived the general  minimax lower bounds for the risk that have not 
been obtained   earlier.  We also  elaborate  on the upper bounds for the risk 
derived earlier by Butucea and Comte (2009).   Using these results, we immediately retrieve  
the upper and the minimax lower  bounds  for the risk of the pointwise  estimator of the mixing density 
with classical and Berkson errors studied by Delaigle (2007). Direct comparison with Delaigle (2007) shows that 
our upper risk bounds are more precise due to more flexible assumptions; they are also   
supplemented by the minimax lower bounds that have not been derived previously.  In addition, our estimator is adaptive.

Furthermore, we expand  our theory in a novel way to incorporate estimation 
of functionals of the form  \fr{eq:Phi} where function $\ph(x)$ does not have a   Fourier 
transform in a regular sense. The new methodology relies on application of inversion formulas.
We show  that by using such formulas, many functionals of interest can be expressed as  
 linear functionals of the real and the imaginary parts of the Fourier transform $\fs$ of $f$ and 
their derivatives. We construct  estimators for such  functionals and derive  the upper bounds for 
their risks. As a particular case of application of our  methodology, we automatically recover     
the estimators of the mixing cumulative distribution function  investigated   by Dattner {\it  et al.} (2011),
as well as their minimax lower and upper bounds for the risk. Other examples include estimation of the $(2M+1)$-th absolute moment 
of deconvolution density and the functional of the form $\iii \te^m (\te^2 +1)^{-1} f(\te) d\te$ with $m \geq 2$. 
As a special case, we obtain  an estimator of the first absolute moment of the deconvolution density
(the problem studied by Cai and Low (2011) in the case of the Gaussian error distribution),
and retrieve   the   upper and lower bounds for the minimax risk of this estimator  as a particular case of our general minimax bounds.

Finally,   we use  the same approach to deal with the situation where deconvolution density $f$ is delta-contaminated
and is defined by \fr{sparseden}. This set-up corresponds to the situation where the vector of unobservable variables  is sparse 
and has, on the average,  only $k_n$ nonzero entries. Following Cai and Low (2011), we assume  that $k_n$  
(or $\mu_n$ in \fr{sparseden}) is known  and the objective is estimating the functional 
over these non-zero entries only.   We show  that convergence rates are determined by the ``effective'' sample size 
$n_\mu = n^{-1} k_n^2$,   propose  an   estimation procedure specifically designed for the sparse case, 
and construct  the  minimax lower and the upper bounds for the risk.
 In particular, an immediate consequence of our theory is a proof that consistent estimation of any linear 
functional of deconvolution density  is impossible if  $k_n = n^{\nu}$ and  $0<\nu \leq 1/2$, since, 
in this case, the effective sample size $n_\mu = n^{2\nu -1} <1$. 
This is a significant  generalization of the results of Cai and Low (2011) who drew the same conclusion specifically for estimation of 
the first absolute moment of the unobservable variable under Gaussian errors.

We also carry out a limited simulation study which compares our estimator $\hPhi_{1,h}$ of the first absolute moment of the mixing density
\fr{abmom} with the estimator $\hPhi_{CL}$ of   Cai and Low (2011). We   show  that, although the two estimators have similar precision, 
 $\hPhi_{1,h}$ is  adaptive to the choice of parameters while $\hPhi_{CL}$ is sensitive to the choice of the parameter $M_{n,0}$,
the maximum of values of $|\te_i|$, $i=1, \cdots, n$. In addition, our simulations confirm that  the precisions of estimation of functionals $\Phi$
and $\Phi_n$ differ  only insignificantly. Indeed, all   results  derived in the paper for estimation of a linear 
functional  $\Phi$ of the deconvolution density can be automatically applied to estimation of the functional
in indirect observations $\Phi_n$.

\section*{Acknowledgments}

Marianna Pensky   was  partially supported by National Science Foundation
(NSF), grants  DMS-1106564 and DMS-1407475. The author  wants to thank Alexandre Tsybakov, Oleg Lepski 
 and  Alexander Goldenshluger  for 
extremely valuable suggestions and discussions.




\newpage



\section{Proofs  }
\label{sec:proofs}
\setcounter{equation}{0}

\subsection{Proof of the adaptive upper bounds for the risk: the standard case}
\label{sec:adap_upper_bounds}

{\bf Proof of Theorem \ref{th:adapt_bounds}. }
Note that the values $\hh_n$ that deliver optimal (or nearly optimal up to a logarithmic factor of $n$) 
are known for all cases except $b= \beta= d = \ga = 0$,    $a <  \af + 1/2$,
can be derived by the straightforward calculus. If  $b= \beta= d = \ga = 0$,    $a <  \af + 1/2$,
denote   
\bes
\psi_h (y) = \frac{1}{2\pi} \int_{-h^{-1}}^{h^{-1}} e^{-i \om y} \frac{\phs (\om)}{\gs(\om)} \, d\om,
\ees
\bes
C_{\psi,1} = \frac{ \Cpht}{\Cgo \pi (\af - a+1)}, \quad
C_{\psi,2} = \frac{ \Cpht^2  \|g\|_\infty}{2\pi (2\af -2a+1)\,\Cgo^2}, \quad
C_{B,0} =   \lkr \frac{2s + 2a-1}{2B^2}\rkr^{\frac{1}{2s+2\af}}. 
\ees
Observe that $\hPhi_h$ can be presented as 
\bes 
\hPhi_h = \frac{1}{n} \sum_{l=1}^n \psi_h (Y_l).
\ees
and also
\be \label{bias_variance}
\| \psi_h \|_\infty \leq \frac{C_{\psi,1}}{h^{(\af - a+1)}}, \quad 
\Var(\hPhi_h) \leq   \frac{C_{\psi,2}}{n\, h^{(2\af -2a+1)}}, \quad
|\EE \hPhi_h - \Phi|  \leq \frac{B  \Cpht \ h^{s + a -1/2}}{\sqrt{2 \pi (2s + 2a-1)}}. 
\ee  
Then, applying  Bernstein inequality to $\sum n^{-1}[\psi_h (Y_l) - \EE \hPhi_h]$  obtain 
\be \label{Bernstein}
\PP \lkr |\hPhi_h - \EE \hPhi_h| > C_{\tau} n^{-1/2} \sqrt{\log n} h^{a-\af -1/2} \rkr 
\leq 2 n^{-\tau} 
\ee 
for any $\tau >0$ and $C_{\tau}$ satisfying the following inequality
\bes 
C_{\tau} \geq \frac{2 \tau C_{\psi,1}}{3} + \sqrt{\lkr\frac{2 \tau C_{\psi,1}}{3}\rkr^2 + 2 \tau C_{\psi,2}}. 
\ees
Set $\tau = 2(2\af -2a +1)$ and 
let $\ho$ and $\jo$ be such that
\bes
\ho = C_{B,0}\, n^{-\frac{1}{2s+ 2\af}}, \quad h_{\jo} \leq \ho < h_{\jo -1}.
\ees
%
Now, note that by definition \fr{hjvalues} of $\hj$, one has
\be \label{Phiersum}
\EE(\hPhhj - \Phi)^2 = \EE [(\hPhhj - \Phi)^2 \II(\hj \leq \jo)] + \EE [(\hPhhj - \Phi)^2 \II(\hj \leq \jo)] \equiv 
\Del_1 + \Del_2.
\ee
Here, 
\beqn \nonumber
\Del_1 &\leq & 2 \EE [(\hPhhj - \hPhi_{h_\jo})^2 \II(\hj \leq \jo)] + 2 \EE [(\hPhi_{h_\jo} - \Phi)^2 \II(\hj \leq \jo)] \\
& \leq & 2  C_\Phi^2 n^{-1}\, \log n\  h_\jo^{-(2\af -2a+1)} + 2 C_{B,0}\, n^{-\frac {2s+2a-1}{2s + 2\af}} 
 \asymp n^{-\frac {2s+2a-1}{2s + 2\af}} \label{Delta1}
\eeqn
In order to find an upper bound for $\Del_2$, denote 
$D_n (h) = 0.5\, C_\Phi \sqrt{\log n/n} h^{a-\af -1/2} - C_{B, \ph} h^{s+a-1/2}$, 
$C_{B, \ph} = B\, \Cpht/(\pi \sqrt{4s + 4a-2})$ and
\bes
\calM_j = \lfi \om:\, \hj = j \rfi = 
\lfi \om:\, \exists \tilj > j\quad {\mbox such\ that} \quad 
|\hPhi_{h_j}- \hPhi_{\htj}| > C_\Phi n^{-1/2} \sqrt{\log n}\ h_{\tj}^{-(\af -a+1/2)} \rfi.
\ees
By formula \fr{bias_variance}, one has
\bes
|\hPhi_{h_j}- \hPhi_{\htj}| \leq |\hPhi_{h_j}- \EE \hPhi_{h_j}| + |\hPhi_{\htj} - \EE \hPhi_{\htj}|
+ C_{B, \ph} \lkr h_j^{s+a-1/2} + h_{\tj}^{s+a-1/2} \rkr.
\ees
Since $\tj > j$, so that $\htj < h_j$, derive
\bes
\PP(\calM_j) \leq \PP (|\hPhi_{h_j}- \EE \hPhi_{h_j}| > D_n (h_j)) + 
\PP (|\hPhi_{\htj}- \EE \hPhi_{\htj}| > D_n (\htj))  
\ees
Moreover, since for $j \geq \jo+1$, one has $h_j \leq \ho/2$, assumption \fr{CPhi} on $C_\Phi$
guarantees that 
\bes
\min \lkr D_n (h_j), D_n (\htj) \rkr \geq C_{\tau} n^{-1/2} \sqrt{\log n}\  h^{a-\af -1/2},
\ees
so that,  by \fr{Bernstein}, one obtains $\PP(\calM_j) \leq 4 n^{-\tau}$. Finally, since one can show by direct
calculations that for $j \leq J$ and some positive absolute constant $C$, one has
\bes
\EE (\hPhi_{h_j} - \Phi)^4 \leq C n^{4\af - 4a} \, (\log n)^{4a-4\af -2},
\ees
derive
\beqn \label{Delta2}
\Del_2 & \leq & \sum_{j = \jo+1}^J \sqrt{\EE (\hPhi_{h_j} - \Phi)^4} \sqrt{\PP(\calM_j)}
\leq C n^{2\af - 2a - \tau/2} \leq C n^{-1}, 
\eeqn
since $\tau = 2(2\af -2a +1)$. 
Combination of formulas \fr{Phiersum}, \fr{Delta1} and \fr{Delta2}  complete the proof.


\subsection{Proof of the lower bounds for the risk: the standard case}
\label{sec:lower_bounds}

{\bf Proof of Theorem \ref{th:lower_bound}. }
The proof is based on application of Theorems 2.1 and 2.2 of Tsybakov (2009)
with the chi-squared divergence.
Introduce function 
\be  \label{Kfungen}
\Ksv (\om) = \lfi 
\begin{array}{ll}
(|\om| -1)^v\, \calP_1 (|\om|), & 1 < \om < 2,   \\
1, & 2 \leq |\om| \leq 3,  \\
(4 - |\om|)^v\, \calP_2 (|\om|), & 3 < \om < 4, \\
0, & \mbox{otherwise}  
\end{array} \right.
\ee
where $\calP_1$ and $\calP_2$ are   such that  $\calP_1 (z) \neq 0$ for $1 \leq z \leq 2$,  
$\Ksv (\om)$ is $(v-1)$ times continuously differentiable  
on the whole real line (i.e., $\calP_1 (2) =  \calP_2 (3) =1$),  and $0 \leq \Ksv (\om) \leq 1$. 
It is easy to see that function $\Ksv (\om)$ is even, real-valued and non-negative.

Consider two functions  
\be \label{f1f2}
f_1 (\te) = b \pi^{-1} (1 + b^2 \te^2)^{-1}, \quad 
 f_2 (\te) =  f_1 (\te) + \lam \fdel (\te)
\ee
where Fourier transform $\fdels (\om)$ of $\fdel (\te)$ is given by 
\be \label{Fourierfd}
\fdels (\om) = \Ks (\om h) \exp \lfi i \arg(\phis (\om))  \rfi
\ee 
and function $\Ks (\om) = \Ks_3 (\om)$ is of the form \fr{Kfungen} with $v = 3$.
Observe that $f_1(\te)$ is a pdf. In order to ensure that $f_2(\te)$ is a pdf, note 
 that, by Assumption on $\rho(\om) = \arg(\phis (\om))$ and  by formula \fr{Kfungen},
function $\fdels (\om)$ is twice continuously differentiable and absolutely integrable if $h < \om_0$, 
hence, $\fdel(\te) \leq \Cdel \te^{-2}$ as $|\te| \to \infty$. 
Therefore, if $\lam$ is small enough, $f_2(\te) \geq 0$. 
Moreover,   $\fdels (0) =  \iii \fdel (\te) d\te =0$, so that $f_2(\te)$ is   also a pdf.

Choose $b$ in $f_1 (\te)$ small enough, so that $f_1 \in \Om_s (B/2) \subset   \Om_s (B)$  
and let 
\be  \label{lamvalue}
\lam^2 = \Cbs^2 h^{2s+1} \quad \mbox{with} \quad \Cbs^2 =   B^2 (s+1) 17^{-s}/4.
\ee
Then, direct calculations yield that $\iii \lam^2 |\fdel(\om)|^2 (\om^2 + 1)^s d\om \leq B^2/4$, and, therefore,
$f_2 \in \Om_s (B)$.

Now, let us evaluate the difference
\beqn \label{difference}
D_h  & = & |\Phi(f_1) - \Phi(f_2)| = \lam \left| \iii \phis (-\om) \fdels (\om) d\om \right| 
= \lam \iii |\phis(\om)| |\Ks (\om h)| d\om.
\eeqn 
Using formula \fr{lamvalue} and Lemma \ref{lem:integral} in Section \ref{sec:suppl} 
 with $A=a/2$, $G = d$, $\aleph =b$, $v=3$ and $l=1$, we derive 
\be \label{Dh_bound}
D_h \geq \lfi
\begin{array}{ll}
C h^{s + a -1/2}, & \mbox{if} \quad d=0,\\
C h^{s + a + 4b -1/2}\, \exp \lkr - d h^{-b} \rkr, & \mbox{if} \quad d>0
\end{array} \right.
\ee

Denote $q_k (y) = \iii g(y -\te) f_k (\te) d\te$,   and let $Q_k(\bY) = \prod_{i=1}^n q_k (Y_i)$,    be
the pdf  of the sample $Y_1 \cdots, Y_n$ under $f_k$, $k=1,2$. We shall find a combination of parameters $n$ and $h$ 
which ensures that the chi-squared divergence $\chi^2(Q_1, Q_2)$ between $Q_1$ and $Q_2$, 
is bounded above. Then,  by Theorems 2.1 and 2.2 of Tsybakov (2009), one has 
$R_n (\Om_s (B)) \geq D_h^2$ where $D_h$ is given by \fr{Dh_bound}.

By Lemma \ref{lem:q1}, one has $q_1(x) \geq \Cbg (x^2 + 1)^{-1}$, and   
Lemma \ref{lem:chi_distance} implies that, in order to ensure   $\chi^2(Q_1, Q_2) \leq \kappa^2$, it is sufficient to establish that  
\be \label{H_integral} 
H = \iii (x^2 + 1)  [q_2(x) - q_1(x)]^2 dx \leq \Cbg \, n^{-1} \log (1 + \kappa^2) 
\ee 
for some $\kappa \in (0,1)$. Note that $H$ in \fr{H_integral} can be written as 
\bes
H = \lam^2 \ \iii (x^2 + 1) [q_\del (x)]^2 dx.  
\ees
Then, $H$ can be split as $H = \lam^2(H_1 + H_2)$ where
\beqns 
H_1 & = & \iii   [q_\del (x)]^2 dx = \frac{1}{2\pi} \iii |\qs_\del (\om)|^2 d\om =  \frac{1}{2\pi} \iii |\gs (\om)|^2 |\Ks(\om h)|^2 d\om
\eeqns
and 
\beqns 
H_2 & = & \iii  x^2\, [q_\del (x)]^2 dx = \frac{1}{2\pi} \iii \left| \frac{d \qs_\del (\om)}{d\om} \right|^2 d\om \leq H_{21} + h^2 H_{22} + H_{23} 
\eeqns
with 
\beqns
H_{21} & = & \frac{1}{2\pi} \ \iii |(\gs)^{'} (\om)|^2 |\Ks(\om h)|^2 d\om \\   
H_{22} & = & \frac{1}{2\pi} \ \iii |\gs  (\om)|^2 |(\Ks)^{'} (\om h)|^2 d\om  
\eeqns
and  $H_{23} = \rho^2 H_1$.
Taking into account condition \fr{gsderiv} and applying Lemma \ref{lem:integral} with 
$A = \af$, $G =  2\ga$, $\aleph = \beta$, $l=2$ and $v=3$ for $H_1$,
with  $A = \af - \tau$, $G =  2\ga$, $\aleph = \beta$, $l=2$ and $v=3$ for $H_{21}$
and $A = \af$, $G =  2\ga$, $\aleph = \beta$, $l=2$ and $v=2$ for $H_{22}$,
obtain that 
\be \label{H_bound}
H \leq \lfi 
\begin{array}{ll}
C \lam^2\, h^{2 \af -1}, & \mbox{if}\quad \ga =0, \\
C \lam^2\, h^{2 \af + U_\tau} \exp \lkr - 2 \ga h^{-\beta} \rkr, & \mbox{if}\quad \ga > 0.
\end{array} \right.
\ee 
Due to the choice of $\lam$ given by  \fr{lamvalue},  combination of formulae \fr{H_integral} and \fr{H_bound} yield 
the following expression for $h = h(n)$:
\bes
h = \lfi 
\begin{array}{ll}
C n^{-\frac{1}{2\af + 2s}}, & \mbox{if}\quad \ga =0, \\
 \lkv (2 \ga)^{-1}\, ( \log n - (2s + 2 \af +  U_\tau+1)\log \log n) \rkv^{-1/\beta}, & \mbox{if}\quad \ga > 0.
\end{array} \right.
\ees
In order to complete the proof, recall that, by Theorems 2.1 and 2.2 of Tsybakov (2009), one has 
$R_n (\Om_s (B)) \geq D_h^2$ where $D_h$ is given by \fr{Dh_bound}.

Finally, consider the case when convergence rates are paarmetric. In this situation, just set 
$\lam = n^{-1/2}$ and  choose $f_1, f_\del \in \Om_s (B/2)$ such that $f_\del$ integrates to zero and 
$f_2 = f_1 + \lam f_\del$ is nonnegative. Then, the lower bounds can be obtained directly from 
Theorems 2.1 and 2.2 of Tsybakov (2009).
 

\subsection{Proofs of  statements  in   Section \ref{sec:challenge}    }

{\bf Proof of Theorem \ref{th:upper_risk_inv}. }
For  derivation of  \fr{risk_inversion}, observe that $\EE \hPhi_1= \Phi_1$ and,
moreover, due to    assumption  A1 and conditions  \fr{cond1_momn} -- \fr{extra_assump}, one has $\Var(\hPhi_1) \leq C n^{-1}$.
For the same reason, the values of $A_{m,j,k,l}$, $l=1,2$, are uniformly bounded, so that $\EE  \widehat{F_m(1)} = F_m(1)$
and $\Var[\widehat{F_m(1)}] \leq C n^{-1}$.
Observe that under assumption \fr{psiup}, one has
\beqns
\Var[\hPhi_{20h}]  & \asymp & n^{-1} \int_1^{h^{-1}} (\om^2 + 1)^{\af - (m+1) a_m} 
\exp \lkr - 2 d \om^b + 2 \ga \om^\beta \rkr \, d\om   \\
 \EE \hPhi_{20h} - \Phi_{20}  & \asymp & \int_{h^{-1}}^\infty (\om^2 + 1)^{- (m+1) a_m}   
\exp \lkr - 2 d \om^b \rkr |\fs(\om)| d\om 
\eeqns
In  order to complete the proof, note that  $\EE \hPhi_h - \Phi = \EE \hPhi_{20h} - \Phi_{20}$ and 
$\Var[\hPhi_h] \leq \Var[\hPhi_{20h}] + C n^{-1}$ and use the Cauchy inequality 
if $f \in \Om_s(B)$  or   upper bounds for $|\fs|$ if $f \in \Xi_s(B)$.
\\


\noindent
{\bf Proof of Theorem \ref{th:absMmom}. }
Observe  that under assumptions of Theorem~\ref{th:absMmom}, 
conditions of Theorem~\ref{th:upper_risk_inv} are satisfied with 
$m= 2M+1$, $a_m =1$ and $b=d=0$. Hence, inequality \fr{risk_inversion} is valid.
To complete the proof of the upper bounds \fr{up_absMmom} for the minimax risk,
use formula \fr{bounds} in Lemma~\ref{lem:rates} with $A_1 = 2M+1+s$, $A_2 = \af - (2M+2)$
and $b=d=0$.


In order to establish the lower bounds \fr{low_absMmom} for the minimax risk,  
use the methodology suggested by Theorems 2.1 and 2.2 of Tsybakov (2009).
Similarly to the proof of Theorem  \ref{th:lower_bound}, choose 
  \be \label{f1f2mom}
f_1 (\te) = \tC\, b   (1 + b^2 \te^2)^{-(2M+3)}, \quad 
 f_2 (\te) =  f_1 (\te) + \lam \fdel (\te),
\ee
where $\fdel (\te) = h^{-1} K_v(\te/h)$,  $\fdels (\om) =   \Ksv (\om h)$ with $v = 4M+4M+7$ 
and $\Ksv(\om)$ is defined by formula \fr{Kfungen}. Note that  function $\fdels (\om)$ is $(4M+6)$ times 
continuously differentiable and absolutely integrable, 
hence, $\fdel(\te) \leq \Cdel \te^{-(4M+6)}$ as $|\te| \to \infty$. 
Moreover,   $\fdels (0) =  \iii \fdel (\te) d\te =0$.
Therefore, if $\lam$ is small enough, $f_2(\te) \geq 0$, so that $f_2(\te)$ is   also a pdf.

Let $B_1$, $B_2$  and $b$ be such that  $f_1 \in \Om_{s} (B_1/2, B_2/2)$.
Note that
\beqns 
\iii \te^{4M+4} f(\te) d\te & \leq &  \iii \te^{4M+4} f_1(\te) d\te + \lam \iii \te^{4M+4} h^{-1} K(\te/h) d\te\\
& \leq & B_1/2 + \lam h^{4M+4} \iii z^4 K(z) dz  \leq B_1
\eeqns
if $\lam$ and $h$ are small enough. Also, by definition of $\Ksv (\om)$, one has 
$|\Ks_v(\om)| \leq 1$ and $|\om|^s |\Ks_v(\om)| \leq 4^s$, thus
\beqns
\lam |\fdels (\om)|\, (|\om|^s +1) = \lam\, |\Ks (\om h)|\, [|\om h|^s h^{-s} + 1] \leq 
\lam\, (4^s h^{-s} +1) \leq B_2/2
\eeqns
provided $\lam = C_{s B} h^s$ for some absolute constant $C_{s B}$ that depends on $s$ and $B_2$ only.
Now, using integration by part and recalling that $\fdels (\om) =   \Ks_v (\om h)$, one derives     
\beqn 
D_h  
& = & |\Phi(f_1) - \Phi(f_2)| =  2 \lam \pi^{-1} (2M+1)!\ \left| \iii   \om^{-(2M+2)}\,  \fdels  (\om) d\om \right| \nonumber\\
& = & 2 \lam \pi^{-1} (2M+1)! h^{2M+1}\,  \int_1^4 z^{-(2M+2)}\, \Ksv (z) dz \geq C \lam h^{2M+1} \geq C\, h^{s+2M+1}.  \label{dist_mom}
\eeqn
Finally, similarly to the proof of Theorem \ref{th:lower_bound}, 
due to Lemma \ref{lem:q1}, one has $q_1(x) \geq \Cbg (x^2 + 1)^{-(2M+3)}$, and   
Lemma \ref{lem:chi_distance} implies that it is sufficient to ensure that 
\be \label{H_mom} 
H = \iii (x^2 + 1)^{2M+3}\  [q_1(x) - q_2(x)]^2 dx = \lam^2 \iii (x^2 +1)^{2M+3} [q_{\del}(x)]^2 dx
 \leq \Cbg \, n^{-1} \log (1 + \kappa^2). 
\ee 
Following the steps in the proof of Theorem \ref{th:lower_bound}, 
we construct an upper bound of $H$ in \fr{H_mom}    as  $H \leq 2^{2M+2} \lam^2(H_1 + H_2)$ with
\beqns 
H_1 & = &  \frac{1}{2\pi} \iii |\gs (\om)|^2 |\Ksv (\om h)|^2 d\om 
\asymp h^{2\af -1 + \beta(2v+1)} \exp \lkr - 2 \ga h^{-\beta} \rkr
\eeqns
and 
\beqn  
H_2 & = &  \frac{1}{2\pi} \iii  \left| \frac{d^{2M+3} \qdels  (\om)}{d\om^{2M+3}} \right|^2 d\om 
\leq \sum_{k=0}^{2M+3} \iii \left| \frac{d^{2M+3-k} \gs  (\om)}{d\om^{2M+3-k}} \right|^2\ 
\left| \frac{d^{k} \Ksv  (\om h)}{d\om^{k}} \right|^2\,  d\om \\
& \asymp & \sum_{k=0}^{2M+3} h^{2k} \iii (\om^2+1)^{-\af + \tau(2M+3-k)} \exp \lkr - 2 \ga |\om|^\beta \rkr
|(\Ks_v)^{(k)} (\om h)| \, d\om. \label{H2expr}
\eeqn 
If $\ga = \beta =0$, then $H_2 \asymp h^{2\af -1}$. Otherwise, apply Lemma \ref{lem:integral} with 
$v =  4M+7 - k$ for the $k$-th term in \fr{H2expr}. Obtain 
$ H_2 \asymp h^{2\af + \beta - 1 + \vartheta} \exp \lkr - 2 \ga h^{-\beta} \rkr$
where $\vartheta = \min( 4M+6 + \beta(4M+8), \beta(4M+15) - \tau(4M+6))$.
Hence,   derive
\be \label{H_bound_mom}
H \leq \lfi 
\begin{array}{ll}
C \lam^2\, h^{2 \af -1}, & \mbox{if}\quad \ga =0, \\
C \lam^2\, h^{2 \af + \beta - 1 + U^*_{\tau}}\, \exp \lkr - 2 \ga h^{-\beta} \rkr, & \mbox{if}\quad \ga > 0.
\end{array} \right.
\ee 
where $U^*_{\tau} = \min(4M+6+\beta(4M+8),\  \beta(8M+15) - \tau(4M+6),\ \beta(8M+14))$.
With  $\lam = C_{s B} h^s$, 
\fr{H_bound_mom} implies that   inequality \fr{H_mom}  holds provided
$h = C n^{-1/(2s+2\af -1)}$ if $\ga=0$, $\beta =0$, and 
$h = [(2\ga)^{-1}( \log n - (2 \af +\beta -1 + U^*_{\tau})\log \log n)]^{-1/\beta}$ if $\ga>0$, $\beta >0$.
In order to complete the proof, recall that, by Theorems 2.1 and 2.2 of Tsybakov (2009), one has 
$R_n (\Om_s (B)) \geq D_h^2$ where $D_h$ is given by \fr{dist_mom}.


\subsection{Proofs of the statements in Section \ref{sec:sparse}    }

{\bf Proof of Theorem \ref{th:sp-low-param}. }
Let $\mun \leq 1/2$. Otherwise, $k_n \geq n/2$ and there is no point considering the case as being sparse. 
Consider two mixing pdfs
\be \label{fmuk}
\fmuk (\te)  = (1 - \mun) \del(\te) + \lam \mun f_k(\te), \quad k=1,2,
\ee 
and corresponding marginal densities
\be \label{qmuk} 
\qmuk (x) =  (1 - \mun) g(x) + \lam \mun q_k(x), \quad k=1,2,
\ee 
where $q_k (y) = \iii g(y -\te) f_k(\te) d\te$  and $\qmuk (y) = \iii g(y -\te) \fmuk (\te) d\te$.
Again, as in the proof of Theorem \ref{th:lower_bound},  let $\Qmuk (\bY) = \prod_{i=1}^n \qmuk (Y_i)$, $k=1,2$,  be
the pdf  of the sample $Y_1 \cdots, Y_n$ under $f_k$, $k=1,2$. Then, by Lemma \ref{lem:chi_distance}, one has
$\chi^2(Q_{\mu, 1}, Q_{\mu, 2}) \leq \kappa^2$  provided 
\be \label{Qmu12}
H_{\mu} = \iii [q_{\mu 1} (x)]^{-1}\ [q_{\mu 1} - q_{\mu 2}]^2\, dx  \leq n^{-1} \log (1 + \kappa^2).
\ee
Direct calculations yield that 
\beqns 
H_{\mu} & \leq & \lam^2 \mun^2 \ \iii \lkv (1 - \mun) g(x) + \lam \mun q_1(x)\rkv^{-1}\ 
\lkv \iii g(x-\te) (f_2(\te) - f_1(\te)) d\te \rkv^2 \ dx\\
& \leq & \frac{\lam^2 \mun^2}{(1 - \mun)^2}\   (I_1 + I_2) \leq 8 C_I\, \lam^2 \mun^2
\eeqns
since $\mun \geq 1/2$ and $(x-y)^2 \leq x^2 + y^2$ for $x,y \geq 0$.
Therefore, $H_{\mu} \leq n^{-1} \log (1 + \kappa^2)$ provided 
$\lam^2 = \lam_0^2 \min(1, (n \mun^2)^{-1})$ where $\lam_0^2 = \log(1 + \kappa^2)/(8 C_I)$. Finally, due to \fr{ineq2}, one has 
\bes
D_h^2 = [\Phi(f_{\mu 1}) - \Phi(f_{\mu 2})]^2 = \lam^2 \Del_{1 2}^2 = \lam_0^2  \Del_{1 2}^2\, \min(1, (n \mun^2)^{-1})
\ees
which, together with Theorems 2.1 and 2.2 of Tsybakov (2009), completes the proof.
\\ 


\medskip

\noindent
{\bf Proof of Corollary \ref{cor:gauss}}.
Before all else, observe that if $\ph(x)$ is constant,  $\ph(x) = C$ for every $x$, 
then $\Phimu=C$ and estimation is unnecessary.

First, consider the case when the even part of $\ph(x)$, $\ph_{even} (x) = 0.5\,(\ph (x) + \ph(-x))$
is not identically equal to zero. Choose two values $\rho_1, \rho_2 \in (0,1)$, $\rho_1 \neq \rho_2$
and set $f_k (\te) = f(\te| \rho_k)$ where $ f(\te| \rho) = \calN (\te| 0, \sig^2 \rho)$, $k=1,2$, the Gaussian pdf  with zero mean and variance  
$\sig^2 \rho_k$. Then, direct calculations yield that $q_k (x) = \calN (x| 0,  \sig^2 (1+\rho_k))$, $k=1,2$.
Hence,  
$$
I_k =  \iii g^{-1}(x) q_k^2 (x)dx = 
\frac{1}{\sqrt{2 \pi  (1 + \rho_k)} \, \sig} \iii \exp \lkr - \frac{x^2 (1 - \rho_k)}{2\sig^2 (1 + \rho_k)} \rkr \, dx ,  k=1,2,
$$
and inequality \fr{ineq1} holds with $C_I = \max \lfi (1-\rho_1)^{-1/2}, (1-\rho_2)^{-1/2} \rfi$. 
Furthermore, note that $A_k = G(\rho_k)$ where
\bes 
G (\rho) = \iii \ph_{even}(\te) f (\te| \rho) d\te 
= \frac{2}{\sig \sqrt{2 \pi \rho}} \ioi z^{-1/2}  \ph  (\sqrt{z})  e^{-\frac{z}{2 \sig^2 \rho}} \, dz.
\ees
Inequality \fr{ineq2} is violated only if $G (\rho)$ takes constant value for $0< \rho < 1$.
It is easy to notice that $\sqrt{\rho}\, G (\rho)$ is proportional to the Laplace transform of 
the function $z^{-1/2}  \ph_1 (\sqrt{z})$. Hence,  $G (\rho)$ is constant if and only if  the Laplace transform of
$z^{-1/2}  \ph_1 (\sqrt{z})$ is equal to $C z^{-1/2}$, which is possible only if $\ph(x)$ is a constant function.
Therefore, if $\ph_{even} (x)$ does not vanish, it is always possible to choose two values, 
$\rho_1 \neq \rho_2$ in $(0,1)$ such that inequality   \fr{ineq2} holds.
 
Now, consider the situation when $\ph_{even} (x) \equiv 0$. Then,  $\ph(\te)$ is an odd function.
Choose $\rho \in (0,1)$ and set $f_1 (\te) =   f(\te| \rho) = \calN (\te| 0, \sig^2 \rho)$. Let   
$f_2 (\te) = f_1(\te)  + f_\del (\te)$ where
\bes
f_\del (\te)=    \mbox{sign\,} (\te)  \, f_1(\te).
\ees
Since, $|f_\del (\te)| \leq f_1(\te)$ for any $\te$, $f_2(\te) \geq 0$. Moreover, 
condition \fr{ineq1} holds for $k=1,2$. Furthermore, since $f_\del (\te)$
integrates to zero, $f_2 (\te)$ is a pdf. It remains to check that 
\bes
\iii \ph (\te) f_\del (\te) d\te = 2 \ioi \ph (\te) f(\te| \rho)  d\te \not\equiv 0 
\ees
for some $\rho \in (0,1)$.
The latter can be accomplished in the same manner as in the case when $\ph_{even} (\te) \not\equiv 0$.
\\

\medskip


\noindent
{\bf Proof of Theorem \ref{th:low_polyn_g}. }
First, let us prove parametric rates. For this purpose, we choose $f_1$ and $f_2$ in 
Theorem \ref{th:sp-low-param} such that  $\fs_k (\om)$ are
$\vso$ times continuously differentiable with   
\be \label{fkconditions} 
\left| \frac{ d^l  \fs_k (\om)}{d \om^l} \right| \leq C_{f}\, |\fs_k (\om)|,\ \ 
l=1,2, \cdots, \vo,  \quad 
\iii (1 + |\om|^{2 \tau \vo})\,  |\fs_k (\om)|^2 d \om < \infty, \quad k=1,2.  
\ee  
Let $\fs_1 (\om)$ and $\fs_2 (\om)$   also satisfy the condition \fr{ineq2}, i.e.
\bes
\Delta_{12} = \frac{1}{2 \pi}\ \iii [\fs_1 (\om) - \fs_2 (\om)] \phis (-\om) d \om \neq 0.
\ees 
It is always possible to find functions $f_1$ and $f_2$   like this.
In order to prove parametric rates \fr{low_bou1}, it is sufficient to show that inequalities \fr{ineq1} 
hold. For this purpose, note that 
$g^{-1} (x) \leq C_{g1}^{-1} (x^2 + 1)^{\vs} \leq C_{\vs}\, C_{g1}^{-1} (x^{2 \vs} + 1)$, 
where constant $C_{\vs}$ depends on $\vs$ only. Hence,
\bes 
I_k = \iii  g^{-1} (x)\  \lkv \iii g(x-\te) f_k(\te) d\te \rkv^2 \,  dx  \leq C_{\vs}\, C_{g1}^{-1} \  (I_{k 1} + I_{k 2}),
\ees
where 
\beqns
I_{k 1} & = & \iii \lkv \iii  g(x-\te) f_k(\te) d\te \rkv^2 \,  dx \leq 
\iii g(x-\te) dx\  \iii g(x-\te) f_k^2 (\te) d\te dx < \infty
\eeqns
and 
\beqns
I_{k 2} & = & \iii  |x|^{2 \vs} \lkv  \iii  g(x-\te) f_k(\te) d\te \rkv^2 \,  dx \leq 
\iii    \lkv  x^{\vso} \iii  g(x-\te) f_k(\te) d\te \rkv^2 \,  dx\\
& = &
\frac{1}{2\pi} \iii \lkv \frac{d^{\vso} [\gs (\om)  \fs_k (\om)] }{d \om^{\vso}}  \rkv^2 d\om 
< \infty		
\eeqns
due to conditions \fr{fkconditions}. Therefore, due to Theorem \ref{th:sp-low-param},
the first three inequalities in \fr{sparse_low_polyng} are valid.
\\

Now, we need to prove nonparametric rates in \fr{sparse_low_polyng}.
 Let, similarly to the proof of Theorem~\ref{th:sp-low-param},  $\fmuk (\te)$ and $\qmuk (x)$ be 
defined in \fr{fmuk} and \fr{qmuk}, respectively,  
where $q_k (y) = \iii g(y -\te) f_k(\te) d\te$  and $\qmuk (y) = \iii g(y -\te) \fmuk (\te) d\te$.
Similarly to the proof of Theorem~\ref{th:lower_bound}, we    consider $f_1 (\te) = b \pi^{-1} (1 + b^2 \te^2)^{-1}$,
so that $f_1 (\te) \in \Om_s (B/4)$ if $b$ is small enough. Then, by Lemma  \ref{lem:q1}, one has $q_1(x) \geq \Cbg (x^2 + 1)^{-1}$.
Let  $ f_2 (\te) =  f_1 (\te) + \lam \fdel (\te)$ where Fourier transform $\fdels (\om)$ of $\fdel (\te)$ is
given by \fr{Fourierfd} with $\Ks (\om) = \Ks_{2 \vo +1} (\om)$ where function $\Ksv (\om)$ is defined by \fr{Kfungen}.  
Again, since $\fdels$ is $(2 \vo)$ times continuously differentiable for $h < \om_0^{-1}$  
and is absolutely integrable, one has $\fdel(\te) \leq \Cdel |\te|^{-2 \vo}$ as $|\te| \to \infty$. 
Therefore, if $\lam$ is small enough, $f_2(\te) \geq 0$ and it is a pdf. 
Also, similarly to the proof of Theorem \ref{th:sp-low-param}, we choose $\lam^2$
of the form  \fr{lamvalue}, so that   $f_2 \in \Om_s (B)$.

Now,   we need to  evaluate the difference $D_h = |\Phi(f_1) - \Phi(f_2)|$ 
given by \fr{difference}.  Using expression for $\lam$ in formula \fr{lamvalue} and Lemma \ref{lem:integral} with 
$A=a/2$, $G = d$, $\aleph =b$, $v= 2 \vo +1$ and $l=1$, we derive 
\be \label{Dh_bound_sp}
D_h \geq \lfi
\begin{array}{ll}
C h^{s + a -1/2}, & \mbox{if} \quad d=0,\\
C h^{s + a + 2b (\vo +1)  -1/2}\, \exp \lkr - d h^{-b} \rkr, & \mbox{if} \quad d>0.
\end{array} \right.
\ee 
It remains to find relation between $n$ and $h$ which ensures that the 
chi-squared divergence $\chi^2(Q_{\mu 1}, Q_{\mu 2})$  is bounded above. 
Lemma \ref{lem:chi_distance} implies that $\chi^2(Q_{\mu 1}, Q_{\mu 2}) \leq \kappa^2$
is guaranteed by 
\be   \label{Hmunew}
H_{\mu} = \lam^2 \mun^2 \iii [g (x)]^{-1}\ (\qdels (x))^2\, dx  \leq n^{-1} \log (1 + \kappa^2).
\ee 
Similarly to the proof of Theorem \ref{th:lower_bound}, we note that 
\bes     
 H_{\mu} \leq C_{\vos} \lam^2 \mun^2 (H_{\vos 1} + H_{\vos 2})
\ees   
where   constant $C_{\vos}$ depends on $\vos$ only and, similarly to the proof of Theorem \ref{th:lower_bound}, 
\beqns 
H_{\vos 1} & = & \frac{1}{2\pi} \iii |\qs(\om)|^2 d\om =  \frac{1}{2\pi} \iii |\gs (\om)|^2 |\Ks(\om h)|^2 d\om
\eeqns
and 
\beqns 
H_{\vos 2} & = & \frac{1}{2\pi} \iii \left| \frac{d^{\vos}  \qs(\om)}{d\om^{\vos}} \right|^2 d\om 
\leq H_{\vos 21} + h^{2 \vos} H_{\vos 22} + H_{\vos 23}
\eeqns
with
\beqns
H_{\vos  21} & = & \frac{1}{2\pi} \ \iii \left| \frac{d^{\vos}  \gs(\om)}{d\om^{\vos}} \right|^2\, |\Ks(\om h)|^2 d\om \\ 
H_{\vos 22} & = & \frac{1}{2\pi} \ \iii |\gs  (\om)|^2 |(\Ks)^{(\vos)} (\om h)|^2 d\om  
\eeqns
and  $H_{\vos 23} = \rho^2 H_{\vos 1}$.

 Take into account condition \fr{gsvsderiv} and apply  Lemma \ref{lem:integral} with 
 $G =  2\ga$, $\aleph = \beta$, $l=2$. Let
$A = \af$  and $v= 2 \vo +1$ for $H_{\vos 1}$  and $H_{\vos 23}$,
$A = \af - \vo \tau$,   and $v=  2 \vo +1$ for $H_{\vos  21}$ and 
$A = \af$  and $v=   \vo +1$ for $H_{\vos  21}$. Then, 
\be \label{H_new_bound}
H_{\mu} \leq \lfi 
\begin{array}{ll}
C \lam^2\, \mun^2\ h^{2 \af -1}, & \mbox{if}\quad \ga =0, \\
C \lam^2\,  \mun^2\  h^{2 \af + U_{\tau, \vo} + 2s +1}\ \exp \lkr - 2 \ga h^{-\beta} \rkr, & \mbox{if}\quad \ga > 0.
\end{array} \right.
\ee 
Due to   \fr{lamvalue},  combination of \fr{Hmunew}, \fr{H_new_bound} and $\mun = n^{\nu -1}$
 yields the following expression for $h = h(n)$:
\bes
h = \lfi 
\begin{array}{ll}
C n^{-\frac{2 \nu -1}{2\af + 2s}}, & \mbox{if}\quad \ga =0, \\
 \lkv \frac{2 \nu -1}{2 \ga}\, \lkr \log n - 
\frac{2 \af + 2s + U_{\tau, \vo} +1 }{2 \nu -1}\ \log \log n\rkr \rkv^{-1/\beta}, & \mbox{if}\quad \ga > 0.
\end{array} \right.
\ees
In order to complete the proof, recall that, by Theorems 2.1 and 2.2 of Tsybakov (2009), one has 
$R_n (\Om_s (B)) \geq D_h^2$ where $D_h$ is given by \fr{Dh_bound_sp}.
\\

\medskip


\noindent
{\bf Proof of Theorem  \ref{th:log_risk}. }
Theorem  \ref{th:log_risk} can be proved by combination of methods used in the proofs of Theorems 
\ref{th:lower_bound} and  \ref{th:low_polyn_g}.  Let $\fmuk (\te)$ and $\qmuk (x)$ be defined in \fr{fmuk} and \fr{qmuk}, respectively,  
where $q_k (y) = \iii g(y -\te) f_k(\te) d\te$  and $\qmuk (y) = \iii g(y -\te) \fmuk (\te) d\te$.
Similarly to the proof of Theorem~\ref{th:lower_bound}, we    consider $f_1 (\te) = b \pi^{-1} (1 + b^2 \te^2)^{-1}$,
so that $f_1 (\te) \in \Om_s (B/4)$ if $b$ is small enough. Then, by Lemma~\ref{lem:q1}, one has $q_1(x) \geq \Cbg (x^2 + 1)^{-1}$.
Let  $ f_2 (\te) =  f_1 (\te) + \lam \fdel (\te)$ where Fourier transform $\fdels (\om)$ of $\fdel (\te)$ is
given by \fr{Fourierfd} with $\Ks (\om) = \Ks_3 (\om)$, where function $\Ksv (\om)$ is given by \fr{Kfungen}.  
Then, both $f_1$ and $f_2$ are pdfs and $D_h = |\Phi(f_1) - \Phi(f_2)| \geq C h^{s + a -1/2}$.
On the other hand, inequality \fr{Qmu12} holds whenever, similarly to \fr{H_integral}, one ensures that 
\bes 
H_{\mu} = \lam^2 \mu_n  \Cbg^{-1}\,  \iii (x^2 + 1)\ (\qdels (x))^2\, dx  \leq n^{-1}   \log (1 + \kappa^2).
\ees
Since $\mun = n^{\nu-1}$, the  latter is guaranteed by choosing 
$$
h =   \lkv (2 \ga)^{-1}\, ( \nu \log n - (2s + 2 \af +  U_\tau+1)\log \log n) \rkv^{-1/\beta}
$$
and yields $D_h \geq C\, (\log n)^{- \frac{2s + 2a -1}{\beta}}$. Application of 
Theorems 2.1 and 2.2 of Tsybakov (2009) completes the proof.


\subsection{Proofs of the supplementary statements  used in Section \ref{sec:lower_bounds}  }
\label{sec:suppl}


\begin{lemma} \label{lem:q1}
Let $\vro \geq 1$ and   $f (\te) =  \Cvro\, b   (1 + b^2 \te^2)^{-\vro}$ where $\Cvro$ depends on $\vro$ only. 
Let $A>0$ be such that  
\bes 
\int_{-A}^A g(x) dx = C_A >0. 
\ees
Then   
\be \label{q_bound}
q  (x) = \iii g(x -\te) f (\te) d\te \geq \Cbg   (x^2 + 1)^{-\vro}
\ee
for some positive constant $\Cbg$ which depends on $\vro$, $\Cvro$,  $b$, $A$ and $C_A$.
\end{lemma}

\noindent
{\bf Proof of Lemma \ref{lem:q1}. }
 Observe that 
\bes 
q(x) \geq \Cvro\, b\ \int_{-A}^A  (1 + b^2(x-z)^2)^{-\vro}\ g(z)\, dz 
\ees
and consider cases $|x|\leq A$ and $|x|> A$ separately.
If $|x|\leq A$, then, due to $|x-z| \leq 2A$ and $(x^2 + 1)^{-1} <1$, one obtains 
\be \label{qbound1}
q(x) \geq  \Cvro b\ (1 + 4 b^2 A^2)^{-\vro}  \int_{-A}^A  g(z)\, dz  \geq \Cvro  C_A b\, (1 + 4 b^2 A^2)^{-\vro} 
\geq \Cvro  C_A  b\, (1 + 4 b^2 A^2)^{-\vro}\,  (x^2 + 1)^{-\vro}.
\ee 
If  $|x| > A$, then, due to $|x-z| \leq 2x$, one derives
\be \label{qbound2}
q(x) \geq  \Cvro b\ (1 + 4 b^2 x^2)^{-\vro}  \int_{-A}^A  g(z)\, dz  \geq \Cvro  C_A  b (1 + 4 b^2)^{-\vro} (1 + x^2)^{-\vro}.
\ee 
In order to obtain \fr{q_bound}, combine \fr{qbound1} and \fr{qbound2} and use 
$\Cbg = \Cvro  C_A b \, \max((1 + 4 b^2 A^2)^{-\vro},\   (1 + 4 b^2)^{-\vro})$.
\\

\medskip


\begin{lemma} \label{lem:integral}
Let function $\Ksv (\om)$ be given by \fr{Kfungen}  
where $\calP_1$ and $\calP_2$ are   such that  $\calP_1 (z) \neq 0$ for $1 \leq z \leq 2$,  
$\Ksv (\om)$ is $(v-1)$ times continuously differentiable  
on the whole real line (i.e., $\calP_1 (2) =  \calP_2 (3) =1$),  and $0 \leq \Ksv (\om) \leq 1$. 
Consider
\be \label{integral}
J_l  = \iii (\om^2 +1)^{-A} \ \exp \lkr -G|\om|^\aleph \rkr |\Ksv (\om h)|^l \, d\om, \quad l=1,2.
\ee
where $A$, $G$ and $\aleph$ are nonnegative constants such that 
$A>0$ and $\aleph =0$ whenever $G=0$. Then, as $h \to 0$, one has
\be \label{asymp_J}
J_l    \asymp \lfi
\begin{array}{ll}
  h^{2A-1}, & \mbox{if}  \quad G=0, \\
   h^{2A + \aleph (l v +1)  -1}\ \exp(-G h^{-\aleph}), & \mbox{if}  \quad G>0. 
\end{array} \right.
\ee
\end{lemma}

\noindent
{\bf Proof of Lemma \ref{lem:integral}. }
It is easy to see that function $\Ksv (\om)$ is even, real-valued and non-negative.
Re-write $J_l$ as 
\bes
J_l = \frac{2}{h} \iof (\om^2 h^{-2} + 1)^{-A}\ \exp \lkr - G h^{-\aleph} \om^\aleph  \rkr\ |\Ksv (\om)|^l \, d\om
\ees
and note that 
\be \label{Jl}
J_l \asymp h^{2A-1}\,  \Uhl \quad \mbox{with} \quad \Uhl =  \iof \exp \lkr - G h^{-\aleph} \om^\aleph \rkr \ |\Ksv (\om)|^l \, d\om.
\ee
If $G=0$, then $\Uhl$ in   \fr{Jl} does not depend on  $h$ and $J_l \asymp h^{2A-1}$.
If $G>0$, the main contribution in \fr{Jl} comes from interval $(1,2)$, so that, introducing new variable $z = \om -1$, obtain 
\bes
\Uhl \approx  \iov  \exp\lkr - G h^{-\aleph} (z+1)^\aleph \rkr\ z^{l v}\, \calP_1 (z+1) dz = 
\iov  \exp\lkr - u_h (z) \rkr\ z^{l v}\, \calP_1 (z+1) dz 
\ees
where function $u_h(z) = G h^{-\aleph} (z+1)^\aleph$ has its minimum on the interval $(0,1)$ at $z=0$.
Note that $u_h (0) =   G h^{-\aleph}$ and $u_h^\prime (z) = \aleph  G h^{-\aleph}\ (z+1)^{\aleph -1}$.
Since $\aleph (z+1)^{\aleph -1}$ lies between $\aleph$ and $\aleph\, 2^{\aleph -1}$,   the mean value theorem implies that
\bes 
\Uhl \asymp \exp (-G h^{-\aleph})\  \iov z^{l v} \calP_1 (z+1) \exp(-  \Caleph G h^{-\aleph}  z)  dz,
\ees
where the constant $\Caleph$  lies between $\aleph$ and $\aleph\, 2^{\aleph -1}$. 
Change   variables $t =  \Caleph  G h^{-\aleph} z$ and recall that $\calP_1$ is a continuous function with $\calP_1 (x) \neq 0$
for $x \in [1,2]$, so that $\calP_1 (x)$ is uniformly bounded above and below for $x \in [1,2]$ and $\calP_1 (z+1) \asymp 1$ for $z \in [0,1]$.
Therefore, 
\be \label{Uhl} 
\Uhl \asymp \exp (-G h^{-\aleph}) h^{l v +1} \  \int_0^{\Caleph  G h^{-\aleph}}    t^{lv}   \exp(- t)  dt \asymp \exp (-G h^{-\aleph}) h^{l v +1}.
\ee 
Combination of \fr{Jl} and \fr{Uhl} complete the proof.
 \\

\medskip


\begin{lemma} \label{lem:chi_distance} 
Let $\bldx = (x_1, \cdots, x_n)$ and let  $Q_k(\bldx) = \prod_{i=1}^n q_k (x_i)$, $k=1,2$,
be two pdfs.  Denote 
\bes    
\chi^2(Q_1, Q_2) = \iii \cdots \iii  Q_1^{-1} (\bldx) \lkv  Q_2 (\bldx) -  Q_1(\bldx) \rkv^2 d\bldx,
\ees
 the chi-squared divergence between $Q_1$ and $Q_2$.
Then, for any $\kappa \in (0,1)$    
\be \label{chi-sufficient}
 \iii [q_1(x)]^{-1}\, [q_1(x) - q_2(x)]^2 dx \leq n^{-1} \log (1 + \kappa^2) \quad 
\mbox{implies} 
  \quad \chi^2(Q_1, Q_2) \leq \kappa^2.
\ee
\end{lemma}

\noindent
{\bf Proof of Lemma \ref{lem:chi_distance}. }
Observe that  
\beqns
I_q & = &  \iii [q_1(x)]^{-1}\, [q_1(x) - q_2(x)]^2 dx = \iii [q_1(x)]^{-1}\, [q_2(x)]^2 dx -1.
\eeqns
Therefore,
\beqns
\chi^2(Q_1, Q_2)  
& = & \iii \cdots \iii Q_1^{-1} (\bldx)\  Q_2^2 (\bldx)  d\bldx - 1 
= \lkv \iii [q_1(x)]^{-1}\, [q_2(x)]^2 dx \rkv^n - 1 \\
& = &
(I_q + 1)^n -1 = \exp[n \log (I_q + 1)] -1.
\eeqns
Taking into account that $\log (1 + z) \leq z$ for $0<z<1$, one obtains
\bes
\chi^2(Q_1, Q_2) \leq \exp(n I_q) -1 \leq \exp[\log(1+ \kappa^2)] -1 = \kappa^2,
\ees
which proves \fr{chi-sufficient}.


\end{document}